\documentclass[12pt]{amsart}
\usepackage{enumerate, amsmath, amsthm, amsfonts, amssymb, mathrsfs, graphicx, paralist, stmaryrd}
\usepackage[all]{xy}
\usepackage[usenames, dvipsnames]{color}
\usepackage[margin=1in]{geometry} 
\usepackage[bookmarks, bookmarksdepth=2, colorlinks=true, linkcolor=blue, citecolor=blue, urlcolor=blue]{hyperref}
\usepackage{enumitem}
\usepackage{tikz-cd}
\usepackage{tikz}
\usepackage{fdsymbol}
\usepackage{mathtools}
\usepackage{caption}

\title[GL-algebras in positive characteristic II]{GL-algebras in positive characteristic II: the polynomial ring}
\author{Karthik Ganapathy}
\address{Department of Mathematics, University of California, San Diego, CA}

\email{\href{mailto:kganapathy@ucsd.edu}{kganapathy@ucsd.edu}}
\urladdr{\url{https://sites.google.com/view/karthik-ganapathy/}}
\usetikzlibrary{%
  matrix,%
  calc,%
  arrows%
}

\newcommand{\cA}{\mathcal{A}}

\newcommand{\cB}{\mathcal{B}}

\newcommand{\bC}{\mathbf{C}}

\newcommand{\rD}{\mathrm{D}}

\newcommand{\bF}{\mathbf{F}}

\newcommand{\bG}{\mathbf{G}}

\newcommand{\bK}{\mathbf{K}}

\newcommand{\bN}{\mathbf{N}}

\newcommand{\bQ}{\mathbf{Q}}

\newcommand{\rR}{\mathrm{R}}

\newcommand{\bS}{\mathbf{S}}
\newcommand{\cS}{\mathcal{S}}
\newcommand{\fS}{\mathfrak{S}}

\newcommand{\bV}{\mathbf{V}}

\newcommand{\bZ}{\mathbf{Z}}

\newcommand{\fa}{\mathfrak{a}}

\newcommand{\fb}{\mathfrak{b}}

\newcommand{\fm}{\mathfrak{m}}

\newcommand{\fp}{\mathfrak{p}}

\renewcommand{\phi}{\varphi}

\newcommand{\lw}{{\textstyle \bigwedge}}

\makeatletter
\def\Ddots{\mathinner{\mkern1mu\raise\p@
\vbox{\kern7\p@\hbox{.}}\mkern2mu
\raise4\p@\hbox{.}\mkern2mu\raise7\p@\hbox{.}\mkern1mu}}
\makeatother

\DeclareMathOperator{\wgt}{wt}
\DeclareMathOperator{\magn}{\text{-mag}}
\newcommand{\pmag}{p\magn}

\DeclareMathOperator{\Div}{Div} 
\DeclareMathOperator{\coker}{coker}

\DeclareMathOperator{\supp}{supp} 
\DeclareMathOperator{\rad}{rad}

\DeclareMathOperator{\cone}{Cone}

\DeclareMathOperator{\ext}{Ext}

\DeclareMathOperator{\End}{End}

\DeclareMathOperator{\Sym}{Sym}

\DeclareMathOperator{\Tor}{Tor}

\DeclareMathOperator{\Spec}{Spec}

\DeclareMathOperator{\Hom}{Hom}

\DeclareMathOperator{\Mod}{Mod}
\newcommand{\id}{\mathrm{id}}

\newcommand{\pol}{\mathrm{pol}}
\newcommand{\Pol}{\mathbf{Pol}}

\newcommand{\tors}{\mathrm{tors}}

\DeclareMathOperator{\Sh}{\bf{\Sigma}}
\DeclareMathOperator{\nSh}{Sh}
\DeclareMathOperator{\De}{\bf{\Delta}}
\DeclareMathOperator{\Rep}{Rep}

\DeclareMathOperator{\Fec}{Vec}

\DeclareMathOperator{\fgen}{fg}
\DeclareMathOperator{\gen}{gen}

\DeclareMathOperator{\chark}{char}
\newcommand{\Sm}{S/\fm^{[q]}}
\newcommand{\Shq}{\Sh_{q/p}}
\newcommand{\Deq}{\De_{q/p}}
\newcommand{\bKq}{{{\bK}_{q/p}}}
\newcommand{\mKq}{\bKq}
\newcommand{\mDeq}{\Deq}
\newcommand{\mShq}{\Shq}

\DeclareMathOperator{\maxdeg}{{maxdeg}}
\DeclareMathOperator{\Dist}{\text{Dist}}
\newcommand{\GL}{\mathbf{GL}}
\newcommand{\GA}{\mathbf{GA}}

\newcommand{\mGaq}{\Gamma_{q/p}}

\newcommand{\fgl}{\mathfrak{gl}}

\makeatletter
\@addtoreset{equation}{section}
\makeatother

\numberwithin{equation}{section}
\newtheorem{theorem}[equation]{Theorem}

\newtheorem{proposition}[equation]{Proposition}
\newtheorem{lemma}[equation]{Lemma}
\newtheorem{corollary}[equation]{Corollary}

\theoremstyle{definition}
\newtheorem{rmk}[equation]{Remark}
\newenvironment{remark}[1][]{\begin{rmk}[#1] \pushQED{\qed}}{\popQED \end{rmk}}
\newtheorem{eg}[equation]{Example}
\newenvironment{example}[1][]{\begin{eg}[#1] \pushQED{\qed}}{\popQED \end{eg}}
\newtheorem{defn}[equation]{Definition}
\newenvironment{definition}[1][]{\begin{defn}[#1]\pushQED{\qed}}{\popQED \end{defn}}

\makeatletter
\renewcommand{\thesubsection}{%
  \ifnum\c@subsection<1 \@arabic\c@section
  \else \thesection.\@arabic\c@subsection
  \fi
}
\makeatother

\theoremstyle{plain}
\newtheorem{mainthm}{Theorem}

\DeclareMathOperator{\slope}{slope}
\date{}
\begin{document}
\begin{abstract}
    We study $\GL$-equivariant modules over the infinite variable polynomial ring $S = k[x_1, x_2, \ldots, x_n, \ldots]$ with $k$ an infinite field of characteristic $p > 0$. We extend many of Sam--Snowden's far-reaching results from characteristic zero to this setting. For example, while the Castelnuovo--Mumford regularity of a finitely generated $\GL$-equivariant $S$-module need not be finite in positive characteristic, we show that the resolution still has finitely many ``linear strands of higher slope". 
   
  The crux of this paper is two technical results. The first is an extension to positive characteristic of Snowden's recent linearization of Draisma's embedding theorem which we use to study the generic category of $S$-modules. The second is a Nagpal-type ``shift theorem" about torsion $S$-modules for which we introduce certain categorifications of the Hasse derivative. 
  These two results together allow us to obtain explicit generators for the derived category. In a follow-up paper, we use the main results here to prove finiteness results for local cohomology modules.
\end{abstract}
\maketitle

\section{Introduction}\label{s:intro}
A $\GL$-algebra over a field $k$ is a commutative $k$-algebra with an action of the infinite general linear group $\GL_{\infty}$ under which it forms a {polynomial} representation.
Around a decade ago, many authors from disparate fields \cite{sno13delta, cef15fi} noticed that $\GL$-algebras exhibit pleasant finiteness properties. Most notably, modules over $\GL$-algebras tend to be noetherian \cite{cef15fi, nss16deg2}, and their free resolutions can often be described using finite data despite being infinite \cite{gl18koszul, ss19gl2}. More recently, Bik--Draisma--Eggermont--Snowden \cite{bdes21geo, bdes23uni} have started systematically studying the spectrum of $\GL$-algebras; Draisma \cite{dra19top} had earlier proved that these spaces satisfy the descending chain condition on $\GL$-stable closed subsets.

In this paper, we analyze the simplest non-trivial $\GL$-algebra $S = k[x_1, x_2, \ldots, x_n, \ldots]$ over an infinite field $k$ of characteristic $p>0$. This complements the project initiated by Sam--Snowden \cite{ss16gl} in which they analyze $S$ over fields of characteristic zero, and our earlier work \cite{gan22ext} on the ``skew" $\GL$-algebra $R = \lw(x_1, x_2, \ldots, x_n, \ldots )$ over $k$. The main point of part I of this series was that the behaviour of $R$ doesn't depend on the characteristic too much, whereas in this paper we'll see that $S$ is much more intricate in characteristic $p$. 

\subsection{Basic structure of the algebra \texorpdfstring{$S$}{S}}
Unlike in characteristic zero, the space of degree $d$ homogeneous polynomials $\Sym^d(\bV)$ is not necessarily irreducible as a $\GL$-representation over $k$. For example, in $\Sym^p(\bV)$, the subspace of $p$-th powers $\langle x_1^p, x_2^p, \ldots, x_n^p, \ldots \rangle$ is a proper $\GL$-subrepresentation. So there are more $\GL$-stable ideals in $S$ over $k$, like the Frobenius powers $\fm^{[q]} \coloneqq (x_1^{q}, x_2^{q}, \ldots, x_n^{q}, \ldots )$ where $q = p^r$.
To make precise how the ideals $\fm^{[q]}$ lead to more complexity in our setting, we recall the notion of a $\GL$-prime ideal in an arbitrary $\GL$-algebra $A$:
    a $\GL$-stable ideal $\fp \subset A$ is \textit{$\GL$-prime} if for every {$\GL$-stable ideal} $\fa, \fb$ with $\fa \fb \subset \fp$, either $\fa \subset \fp$ or $\fb \subset \fp$, i.e., $\GL$-prime ideals are the prime ideals of $A$ when considered as a commutative algebra object in $\Rep^{\pol}(\GL)$. In Section~\ref{s:basicresults} we show:
\begin{mainthm}[$\GL$-spectrum of $S$]\label{thm:introspec}
    Over $k$, the $\GL$-prime ideals of $S$ are $(0)$ and $\fm^{[p^r]}$ with $r \in \bN$.
\end{mainthm}
In contrast, the $\GL$-prime ideals of $S$ in characteristic zero are just the ones we expect, namely the maximal ideal $\fm$ and the zero ideal.  While the above theorem is not too difficult to prove, describing the $\GL$-prime ideals of an arbitrary $\GL$-algebra is quite challenging: Snowden \cite{sno20spe} showed that supergeometry explains everything in characteristic zero, but in characteristic $p$, one must somehow incorporate Frobenius as well.

\subsection{\texorpdfstring{$\GL$}{GL}-equivariant \texorpdfstring{$S$}{S}-modules}
For us, an $S$-module will always mean a \textit{$\GL$-equivariant $S$-module}; an $S$-module is \textit{finitely generated} if it is generated (as an ordinary $S$-module) by the $\GL$-orbit of finitely many elements. A key result which led to the explosion of interest in $\GL$-algebras is a local noetherianity result, which also holds for $S$ in characteristic $p$ (Theorem~\ref{thm:locnoeth}). 

The main structural result of our paper, in light of Theorem~\ref{thm:introspec}, is the existence of a ``derived prime cyclic filtration" for $S$-modules. We let $\rD^b_{\fgen}(\Mod_S)$ be the bounded derived category of finitely generated $S$-modules. For a partition $\lambda$, let $L_{\lambda}$ be the irreducible $\GL$-representation with highest weight $\lambda$.
\begin{mainthm}[Generators for $\rD^b_{\fgen}(\Mod_S)$]\label{thm:gensdbmodintro}
       The bounded derived category $\rD^b_{\fgen}(\Mod_S)$ is generated, as a triangulated category, by the classes of modules $S \otimes L_{\lambda}$ and $S/\fm^{[p^e]} \otimes L_{\lambda}$ as $\lambda$ varies over all partitions and $e$ varies over all integers.       
\end{mainthm}
We prove this in Section~\ref{s:structure} of the paper after proving a similar result for $\Sm$ (Proposition~\ref{prop:gensdbSm}). Much of this paper is dedicated to proving Theorem~\ref{thm:gensdbmodintro}, so we sketch the strategy of proof in Section~\ref{ss:strategy}. While this theorem is admittedly somewhat technical, it is a very powerful result. Indeed, we can swiftly deduce novel finiteness results for the resolutions of $S$-modules as its corollaries.

\subsection{Resolutions of \texorpdfstring{$S$}{S}-modules}
Over fields of characteristic zero, Sam--Snowden showed that the minimal free resolution of an $S$-module exhibits strong finiteness properties. For instance, finitely generated $S$-modules have finite Castelnuovo--Mumford regularity \cite[Corollary~6.3.5]{ss16gl}, i.e., there exists an integer $\rho(M)$ such that the $i$-th syzygy of $M$ is concentrated in degree $\leq i + \rho(M)$ for all $i \geq 0$. This readily fails in characteristic $p$, as the next example shows.
\begin{example}\label{exm:koszul}
 The ideal $\fm^{[p]}$ is generated by a regular sequence, so the Koszul complex
\[ \ldots \to S \otimes \lw^i(\bV^{(1)}) \to \ldots \to S \otimes \lw^2(\bV^{(1)}) \to S \otimes \bV^{(1)} \to S \to 0\]
is a resolution of $S/\fm^{[p]}$; here $\bV^{(1)} = \langle x_1^p, x_2^p, \ldots, x_n^p, \ldots \rangle$ is the ``Frobenius twist" of the standard $\GL$-representation $\bV \coloneqq k^{\infty}$. In particular, $\Tor_i^S(S/\fm^{[p]}, k)$ is concentrated in degree $pi$, so $S/\fm^{[p]}$ has infinite regularity. 
\end{example}
Before stating our result, we note that an $S$-module has finite regularity if and only if all nonzero entries in its \textit{Betti table} are concentrated in finitely many rows (= line of slope $0$), so our result generalizes Sam--Snowden's theorem.
\begin{mainthm}\label{thm:bettiS}
    The Betti table of a finitely generated $S$-module is supported in finitely many lines (of possibly nonzero slope), i.e., there exist finitely many lines in the Betti table of $M$ such that all entries not on  these lines are zero.
\end{mainthm}
\begin{proof}[Proof of Theorem~\ref{thm:bettiS} assuming Theorem~\ref{thm:gensdbmodintro}]
Similar to the argument in Example~\ref{exm:koszul}, the module $S/\fm^{[q]}$ is resolved by the Koszul complex. Using this observation, we see that the explicit generators in Theorem~\ref{thm:gensdbmodintro} have Betti tables concentrated in one line: the Betti table of $S/\fm^{[q]} \otimes L_{\lambda}$ is concentrated in one line of slope $q-1$ as the $(i, j)$-th entry is nonzero if and only if $j = |\lambda| + (q-1)i$. Theorem~\ref{thm:bettiS} now follows by standard d{\'e}vissage arguments.
\end{proof}
Theorem~\ref{thm:bettiS} suggests that there are more hidden structures on $S$-modules. We pursue this in a follow-up paper \cite{gan24lnnr}.
\begin{figure}\label{fig:bettitable}
\begin{displaymath}
\begin{matrix}
   & 0 & 1 & 2 & 3 & 4 & 5 & \cdots \\
0: & \spadesuit & . & . & . & . & . &   \\
1: & . & . & . & . & . & . & \\
2: & . & \spadesuit & . & . & . & . & \\
3: & . & . & . & . & . & . & \\
4: & . & . & \spadesuit & . & . & . & \\
5: & . & . & . & . & . & . & \\
6: & . & . & . & {\spadesuit}  & . & . & \\
\vdots & . & . & . & . & . & . & . &
\end{matrix}
\hspace{1cm}
\begin{matrix}
   & 0 & 1 & 2 & 3 & 4 & 5 & \cdots\\
3: & {\vardiamondsuit} & {\vardiamondsuit} & {\vardiamondsuit} & {\vardiamondsuit} & {\vardiamondsuit} & {\vardiamondsuit} & {\vardiamondsuit} \\
4: & . & . & {\spadesuit} & . & . & . &   \\
5: & . & . & . & {\spadesuit} & . & . & \ldots\\
6: & . & . & . & . &{\spadesuit}  & . &  \\
7: & . & . & . & . & . & {\spadesuit} &   \\
\vdots & . & . & . & . & . & . & \ddots 
\end{matrix}
\end{displaymath}
\caption*{(Illustration of Theorem~\ref{thm:bettiS})
All entries not represented by the suits vanish in the Betti tables of $S/\fm^{[3]}$ with $p = 3$ (left) and of $\fm^{[p]}\fm$ with $p = 2$ (right); see \cite{murai20betti} for similar results.
}
\end{figure}

\subsection{Strategy and technical results}\label{ss:strategy}
We now explain the strategy we employ to obtain Theorem~\ref{thm:gensdbmodintro}. For a $\GL$-prime ideal $\fp \subset S$, we first study $S$-modules that are supported at $\fp$. There are essentially two distinct cases: when $\fp = (0)$ and $\fp = \fm^{[q]}$ for $q = p^r$. The results of Section~\ref{s:snowdenshift} handle the former case, and the results of Section~\ref{s:nagpalshift} handle the latter. 
We begin with an overview of these two sections, highlighting how the results in positive characteristic, though similar to those in characteristic zero, present significantly greater technical challenges---and how we address them.

\subsubsection{Embedding theorem}\label{ss:embedintro}
To study torsion-free $S$-modules (i.e., modules supported at the zero ideal), we prove an ``embedding" theorem for finitely generated \textit{polynomial} $\GL$-algebras, i.e., $\GL$-algebras of the form $\Sym(W)$; recall that $S$ is the polynomial $\GL$-algebra with $W$ the standard representation $\bV$. Snowden \cite{sno21stable} proved the same result in characteristic zero to study the ``fraction fields" of these algebras; we extend the embedding theorem to all characteristics. 
\begin{mainthm}[Embedding theorem]\label{thm:embeddingintro}
   Let $W$ be a finite length $\GL$-representation and let $M$ be a torsion-free $\Sym(W)$-module. The module $M$ embeds into a finitely generated flat module, i.e., there exists a finitely generated flat $\Sym(W)$-module $F$ along with an injection $M \to F$. 
\end{mainthm}
Snowden's proof crucially relies on a simple observation about $\GL$-representations in characteristic zero: every representation has a nonzero weight vector of weight $1^n$, and so any subrepresentation of a tensor product $V \otimes W$ contains a \textit{tensor-disjoint} weight vector, i.e., a weight vector of the form $\sum v_i \otimes w_i$ with the weight of $v_i$ having disjoint support from that of $w_i$ for all $i$. However, $\GL$-representations over $k$ need not have weight vectors of weight $1^n$ (see also Remark~\ref{rmk:bkprest}).
Our main technical lemma (Proposition~\ref{prop:disjointweights}) is the existence of tensor-disjoint weight vectors in characteristic $p$ using 
a vanishing theorem of Friedlander--Suslin \cite[Theorem~2.13]{fs97coh} for the Ext groups between additive polynomial functors and nontrivial tensor products. 
see Section~\ref{ss:flatwgts} for precise results. 

The rest of Section~\ref{s:snowdenshift} is quite similar to Snowden's work \cite{sno21stable} with one major exception stemming from $\Rep^{\pol}(\GL)$ not being semisimple over $k$. We get around this by exploiting the noetherianity of certain Rees algebras and the injectivity of $\Sym(\bV)$ in $\Rep^{\pol}(\GL)$.
\subsubsection{Shift theorem}\label{sss:shiftintro}
In Section~\ref{s:nagpalshift} we shift gears and study $\Sm$-modules since modules with support $\overline{\fm^{[q]}} \subset \Spec_{\GL}(S)$ can be filtered by $\Sm$-modules.
This section is devoted to proving a ``shift theorem" for $\Sm$-modules inspired by the beautiful theorem of Nagpal \cite{nag15fi}. Similar to our earlier paper on the exterior algebra \cite{gan22ext}, our argument is modelled on Li--Yu's \cite{ly17fi} simplification, although the proof is considerably more complicated for $\Sm$. 

We define a sequence of exact endofunctors $\{\Sh_n\}_{n \in \bN}$ of $\Mod_{\Sm}$ which generalizes the classical Schur derivative functor $\Sh \coloneqq \Sh_1$. The sequence $\{\Sh_n\}$ is related to the Schur derivative $\Sh$ in the same way that the Hasse derivatives (of polynomials) are related to the ordinary derivative (see Sections~\ref{ss:hasseschur} and~\ref{ss:shiftdefn}). We then prove:
\begin{mainthm}[Shift theorem for $\Sm$]\label{thm:shiftintro}
   Let $M$ be a finitely generated $\Sm$-module. For $l \gg 0$, the $\Sm$-module $\mShq^l(M)$ is flat.
\end{mainthm}
Forgetting the $\GL$-action, a flat $\Sm$-module is a free $\Sm$-module, so its structure is not too complicated. We sketch the proof of the shift theorem here to illustrate the difficulties we encounter, which may otherwise be obscured in Section~\ref{s:nagpalshift}.

Given a finitely generated $\Sm$-module $M$ it is easy to see that $\mShq^l(M)$ is torsion-free for $l\gg 0$, so we may assume from the outset that $M$ is torsion-free. In this case, let
\begin{displaymath}
    \to F_i \to F_{i-1} \to \ldots \to F_2 \to F_1 \to F_0 \to M \to 0
\end{displaymath}
be a \textit{degree-minimal} flat resolution of $M$; a degree-minimal flat resolution is the closest approximation to a minimal projective resolution for $\GL$-algebras over $k$ (see Example~\ref{exm:koszul2}). Applying the functor $\mDeq \coloneqq \coker(\id \to \mShq)$ to the above resolution results in a flat resolution of $\mDeq(M)$:
exactness is preserved since $M$ is torsion-free, and the functor $\mDeq$ takes flat modules to flat modules. 
This permits us to deduce results about (the resolution of) $M$ using facts about $\mDeq(M)$ which is less complicated being generated in lower degree. In many other contexts, cokernel functors like $\mDeq$ also preserve degree-minimality, a critical property used by Li--Yu \cite{ly17fi}. This property fails for the action of $\mDeq$ on $\Sm$-modules since $\mDeq$ annihilates modules generated by $r$-fold Frobenius twisted representations (Corollary~\ref{cor:dfrobzero}).

There are two ideas we use to overcome these issues. First, we classify the finitely generated torsion-free modules for which $\mDeq$ vanishes in Section~\ref{ss:deltavanishing}. This idea is essentially the same as what appears in our earlier paper \cite{gan22ext} albeit harder to execute (for example, compare Lemma~\ref{lem:dmzero} with \cite[Lemma~4.9]{gan22ext}). The second idea is about controlling the numerical properties of the degree-minimal flat resolution of $\mShq^l(M)$ for $l \gg 0$, which we believe is a novel strategy.
By induction on the generation degree, we can assume that $\mDeq(M)$ satisfies the conclusions of the shift theorem. So by replacing $M$ with $\mShq^l(M)$, we can also assume that $\mDeq(M)$ is flat. Despite the fact that $\mDeq$ does not preserve degree-minimality, we can still control the \textit{slope} of $M$, defined as
\[ \slope(M) \coloneqq \sup_i \frac{t_i(M) - t_0(M)}{i},\]
where $t_i(M) \coloneqq \maxdeg \Tor_i^{\Sm}(M, k)$, provided that the resolution of $M$ is not degenerate to begin with (this is the first case in the proof of Lemma~\ref{lem:t1lesst0}). Specifically, we show that if $\mDeq(M)$ is flat, then $\slope(M) \le q/p - 1$ (the second case in the proof of Lemma~\ref{lem:t1lesst0}). However $\Sm$ is a complete intersection of forms of the same degree $q$, so by Proposition~\ref{prop:cibound}
if $\slope(M) < q/2$ then $M$ must be flat! By induction, we obtain Theorem~\ref{thm:shiftintro}. This idea is executed in Section~\ref{ss:pfshift}. 
\subsubsection{Proving Theorem~\ref{thm:gensdbmodintro}}
We finally explain how to prove Theorem~\ref{thm:gensdbmodintro}. Given an $S$-module $M$, using the results about the generic category in Section~\ref{ss:genericSmod}, we obtain a complex
\[ 0 \to M \to F_0 \to F_1 \to \ldots \to F_r \to 0\]
such that each $F_i$ is a finitely generated flat module and the homology of the complex is torsion (Proposition~\ref{prop:resolutionthmS}). We also prove the existence of a complex of $\Sm$-modules satisfying similar properties (Proposition~\ref{prop:resolutionSm}), but this follows by (essentially) iterating Theorem~\ref{thm:shiftintro}. These two results along with some general nonsense (Lemma~\ref{lem:derivedlemma})
give Theorem~\ref{thm:gensdbmodintro}; see Section~\ref{ss:proofBE} for details. 
\subsection{Relation to other work}
\subsubsection{Sam--Snowden's work in characteristic zero}
Our paper can be seen as the positive characteristic analogue of Sam--Snowden's influential paper \cite{ss16gl}. However, many of their results have to be adjusted appropriately, as we have already seen with Theorem~\ref{thm:bettiS}. For example, projective $S$-modules only have finite injective dimension (over $\bQ$, they are injective on the nose); the Krull--Gabriel dimension of $\Mod_S$ is infinite (it is one-dimensional over $\bQ$); and finitely generated $S$-modules need not have finite injective dimension (over $\bQ$, they do). These create fundamental obstacles to generalize some of their deeper results like Koszul duality for $S$-modules; see \cite[Section~4]{gan24lnnr}.

Sam--Snowden perform a similar analysis for the $\GL$-algebras $\Sym(\bV^{\oplus n})$ in \cite{ss16gl}. We have not pursued the positive characteristic story, as merely classifying their $\GL$-prime ideals seems difficult (see Remark~\ref{rmk:glprimesminors}).

\subsubsection{The infinite general affine group}
The generic category $\Mod_S^{\gen}$ is equivalent to the category of polynomial representations of the infinite general affine group $\GA$. We have not included a proof of this here, but it is essentially similar to the characteristic zero proof \cite{ss19gl2}. This viewpoint allows us to investigate questions of a representation theoretic flavor that are otherwise not apparent. 
Recently, Kriz \cite{kriga24} defined interpolation categories $\Rep(\GA_t(\bC))$ (in the sense of Deligne) for general affine groups using representation stability. It would be interesting to see if such constructions can also be extended to positive characteristic as well.
\subsubsection{$\GL$-varieties}
Bik--Draisma--Eggermont--Snowden have initiated the systematic study of $\GL$-varieties \cite{bdes21geo, bdes23uni} which are the (ordinary) spectrum of $\GL$-algebras. They have proved remarkable technical results about $\GL$-varieties over $\bC$ and a subset of their group also generalized some results to characteristic $p$ \cite{bds24charp}. The theory of $\GL$-equivariant modules over an arbitrary $\GL$-algebra will be the correct notion of a ``quasi-coherent sheaf" on the corresponding $\GL$-variety; we must stress that the geometric theory developed so far does not use any scheme-theoretic notions as of now.
\subsubsection{Resolutions of symmetric ideals}
A central goal in commutative algebra now is to analyze the asymptotic behaviour of invariants of ``compatible" sequences of modules $\{M_n\}$ over $k[x_1, x_2, \ldots, x_n]$.
Along these lines, Le--Nagel--Nguyen--R{\"o}mer have a conjecture stating that the regularity and projective dimension should grow linearly in $n$ when $M_n$ is an $\fS_n$-equivariant ideal in $k[x_1, x_2, \ldots, x_n]$ with $M_n \subset M_{n+1}$ \cite{lnnr20pdim, lnnr21reg}. The conjecture is known in special cases: for $\fS_n$-equivariant \textit{monomial} ideals \cite{murai20betti, rai21reg} and $\GL_n$-equivariant \textit{modules} when $\chark(k) = 0$ \cite{ss16gl, ss22pdim}. We extend the latter result to positive characteristic in \cite{gan24lnnr} and believe our approach can be generalized to a larger class of ideals.
\subsubsection{Stable cohomology} 
The cohomology of line bundles on flag varieties exhibits peculiar phenomena in characteristic $p$ when compared to the characteristic zero side, which is controlled by the Borel--Weil--Bott theorem. Raicu--VandeBogert \cite{rv23flag} made a breakthrough recently when they recognized that for a fixed $i$, the $i$-th cohomology of certain sequences of line bundles $\mathcal{O}(\lambda^{[n]})$ on $Fl(k^n)$ stabilizes to a polynomial functor as $n \to \infty$. While we have not made an explicit connection to their theory, we believe that at least some techniques from this paper will be pertinent. By way of an early example, along with Raicu, we resolve Conjecture~6.4 from \cite{grv23opac} for which the Hasse--Schur derivatives from Section~\ref{ss:hasseschur} enable us to make key reductions. The conjecture concerns the irreducible $\GL$-representations that occur in degrees $\geq n+p$ of $B \otimes B_n$ where  $B \coloneqq S/\fm^{[p]}$. It would be particularly interesting if this conjecture can be related to structural results about $\Mod_B^{\gen}$  like block theory.

\subsection{Notation}
\begin{description}[align=right,labelwidth=2.5cm,leftmargin=!]
\item[$k$] an infinite field of characteristic $p > 0$
\item[$\bV$] the infinite dimensional $k$-vector space with basis $\{e_i\}_{i \ge 1}$
\item[$\GL$] the group of automorphisms of $\bV$ fixing all but finitely many of the $e_i$
\item[$\Rep^\pol(\GL)$] the category of polynomial representations of $\GL$
\item[$L_{\lambda}$] the irreducible polynomial representation of $\GL$ with highest weight $\lambda$
\item[$-^{<n}$] the submodule generated by all elements of degree less than $n$
\item[$-^{(s)}$] the $s$-fold Frobenius twist of a $\GL$-representation/polynomial functor
\item[$t_i(-)$] the generation degree of $\Tor_i^A(-, k)$ (the algebra $A$ is suppressed in this notation as it is clear in context)
\end{description}

\subsection*{Acknowledgements}{The author is grateful to Andrew Snowden for being generous with his ideas and to Teresa Yu for a close reading of this paper. The author was supported by the University of Michigan Prasad Family Fund and NSF grants DMS \#1453893, DMS \#1840234, DMS \#2301871.}

\section{Preliminaries}\label{s:preliminary}
We refer the reader to \cite[Section~2]{gan22ext} for details on polynomial representations of $\GL$ and strict polynomial functors relevant to the contents of this paper. 
We emphasize a nonstandard notation we use starting in this section: given a polynomial functor $F$ and a $k$-vector space $V$, we let $F\{V\}$ be the result of applying $F$ to $V$. Similarly, given a $\GL$-representation $M$ and $k$-vector space $V$, we let $M\{V\}$ be the result of applying $M$, regarded as a polynomial functor, on $V$. 

Throughout this section, $A$ will denote an arbitrary $\GL$-algebra with $A_0=k$, and $S$ will be the $\GL$-algebra $\Sym\{\bV\}$. An $A$-module refers to a module object for $A$ in $\Rep^{\pol}(\GL)$. Spelled out, an $A$-module $M$ is the data of a $k$-vector space $M$ with an action of the $k$-algebra $A$ and a linear action of $\GL$ such that $M$ is a polynomial representation of $\GL$ and the two actions satisfy $g(am) = g(a) g(m)$ for all $g \in \GL, a \in A$ and $m \in M$.

\subsection{Actions on a module} \label{ss:actions}
Given an $A$-module $M$ and a vector space $V$, the $\GL(V)$-representation $M\{V\}$ is evidently a $\GL(V)$-equivariant module over $A\{V\}$. Since these representations are polynomial, the action of $\GL(V)$ extends to an action of $\End(V)$ on $A\{V\}$ and $M\{V\}$. Assume now that $V$ is a finite-dimensional vector space. In this case, the action of $\GL(V)$ on $M\{V\}$ also induces an action of the distribution algebra $\Dist(\GL(V))$ on $M\{V\}$ (see \cite[Chapter~7]{jan03alg} for the definition of the distribution algebra). All these actions evidently extend to $M$.

In the case when $V = k^n \subset \bV$, we provide a concrete description taken from \cite[Section~2.2]{kuj06stpt} of $\Dist(\GL_n)$ arising from the Kostant $\bZ$-form on the universal enveloping algebra of $\GL_n$.
Let $\fgl_n(\bC)$ be the lie algebra of $\GL_n(\bC)$ with $e_{i,j}$ being the $(i,j)$-th matrix unit, and let $U(\GL_n(\bC))$ be the universal enveloping algebra. Define $U(\GL_n)_{\bZ}$ as the $\bZ$-subalgebra of $U(\GL_n(\bC))$ 
generated by the elements $e_{i,j}^{(l)} \coloneqq \frac{e_{i,j}^l}{l!}$ and 
$\binom{h_i}{l}\coloneqq
\frac{e_{i,i}(e_{i,i}-1)\ldots(e_{i,i}-l+1)}{l!}$ (where $i, j, l \in \bN$ with $i\ne j$, $1 \le i,j\le l$, and $l \geq 1$). We have an isomorphism of Hopf algebras,
\begin{displaymath}
    \Dist(\GL_n) \cong k \otimes_{\bZ} U(\GL_n)_{\bZ} 
\end{displaymath}
and we abuse notation by using the symbols $e_{i,j}^{(l)}$ and $\binom{h_i}{l}$ for the images of these elements of $U(\GL_n)_{\bZ}$ in $\Dist(\GL_n)$ under the above isomorphism.

\begin{remark}
    For $s, l \in \bN$, we have $e_{i,j}^{(l)} e_{i,j}^{(s)} = \binom{s+l}{l} e_{i,j}^{(l+s)}$, and similarly, $(l+1) \binom{h_i}{l+1} = \binom{h_i}{l}({h_i} - l) $. Therefore, $\Dist(\GL_n)$ is generated by the elements $e_{i,j}^{(p^r)}$ and $\binom{h_i}{p^r}$ with $1 \leq i,j \leq n$, $i \ne j$ and $r \in \bN$.
\end{remark}

\begin{example}\label{exm:comodule}
Let $\Delta$ be the comodule map $\Dist(\GL_n) \to \Dist(\GL_n) \otimes \Dist(\GL_n)$. In $U(\GL_n)_{\bZ}$, we have $\Delta(e_{i,j}^{(l)}) = \frac{\Delta(e_{i,j}^l)}{l!} = \frac{(e_{i,j} \otimes 1 + 1 \otimes e_{i,j})^l}{l!}$. So reducing mod $p$, we get 
\[\Delta(e_{i,j}^{(l)}) = \sum_{s=0}^l e_{i,j}^{(s)} \otimes e_{i,j}^{(l-s)}\]
in $\Dist(\GL_n)$.
\end{example}
\begin{lemma}\label{lem:hassevanishing}
Assume $A$ is a $\GL$-algebra over $k$. Let $x$ be a homogeneous element of positive degree in $A$ and $q$ be a power of $p$. For arbitrary $i, j, n$ with $q \nmid n$, we have $e_{i, j}^{(n)}(x^q) = 0$.
\end{lemma}
\begin{proof}
    The proof of \cite[Lemma~2.9]{ess19big} applies.
\end{proof}

\subsection{Flat modules and resolutions} \label{ss:flatres}
We recall some structural results about flat modules over $\GL$-algebras \cite[Section~2.4]{gan22ext} and introduce \textit{degree-minimal} flat resolutions for $A$-modules. We first recall a crucial definition:
\begin{definition}
An $A$-module is \textit{induced} if it is isomorphic to $A \otimes W$ for some polynomial representation $W$. An $A$-module is \textit{semi-induced} if it has a finite filtration where the successive quotients are induced modules.
\end{definition}

Given an $A$-module $M$, let $\Tor_i^A(M, -)$ be the left derived functors of the right exact functor $M \otimes_A -$. 
An $A$-module $M$ is \textit{flat} if the functor $M \otimes_A -$ is exact. 
We let $t_i(M) = \max\deg \Tor_i^A(M, k)$ as a polynomial representation with the convention that $\max\deg(0) = -1$; we suppress the algebra $A$ in this notation as it is usually clear which $\GL$-algebra we are working with.

By Nakayama's lemma, the module $M$ is generated in degrees $\le n$ if and only if $t_0(M) \le n$. We denote by $M^{< n}$ the $A$-submodule of $M$ generated by all elements of degree $< n$. If $t_0(M) = n$, the quotient $M/M^{<n}$ will be generated in degree $n$; we implicitly use this throughout the paper.

The proofs of the next three results are in \cite{gan22ext}.
\begin{proposition}\label{prop:flatequalssemi}
    Let $M$ be a finitely generated $A$-module. The following are equivalent:
    \begin{enumerate}[label=(\alph*)]
        \item $M$ is semi-induced,
        \item $M$ is flat,
        \item $\Tor_i^A(M, k) = 0$ for all $i > 0$, and
        \item $\Tor_1^A(M, k) = 0$.
    \end{enumerate}
\end{proposition}

\begin{lemma}\label{lem:relationsinlowdeg}
     Let $M$ be an $A$-module generated in degree $n$ such that $t_1(M) \leq n$. The natural map $A \otimes M_n \to M$ is an isomorphism.
\end{lemma}

\begin{corollary}\label{cor:semises}
Let $0 \to M_1 \to M_2 \to M_3 \to 0$ be a short exact sequence of finitely generated $A$-modules.
\begin{enumerate}[label=(\alph*)]
    \item If $M_1$ and $M_3$ are semi-induced, then so is $M_2$.
    \item If $M_2$ and $M_3$ are semi-induced, then so is $M_1$.
\end{enumerate}
\end{corollary}

In positive characteristic, minimal projective resolutions need not always exist since the underlying category $\Rep^{\pol}(\GL)$ is not semisimple. We work with a certain weakening of this notion in this paper. The next definition appears in \cite{ly17fi} as ``adaptable" resolution.
\begin{definition}
A flat resolution $\cdots \rightarrow F_i \rightarrow F_{i-1} \rightarrow \cdots \to F_0 \rightarrow M = F_{-1} \rightarrow 0$ of an $A$-module $M$ is \textit{degree-minimal} if $t_0(F_i) = t_0(\ker(F_{i-1} \rightarrow F_{i-2}))$ for all $i \geq 1$. 
\end{definition}
The usual argument shows that finite degree-minimal flat resolutions exist for locally noetherian $A$.
\begin{lemma}
Assume $\Mod_A$ is locally noetherian. Let $M$ be a finitely generated $A$-module. Then $M$ has a degree-minimal flat resolution by finitely generated flat $A$-modules.
\end{lemma}
\begin{proof}
    We have a surjective map $F_0 \coloneqq \bigoplus_{i=0}^{t_0(M)} A \otimes  M_i \to M$. The module $F_0$ is induced, thus flat and finitely generated, since each $M_i$ is a finite length $\GL$-representation for all $i$. The kernel of this map is finitely generated as $\Mod_A$ is locally noetherian and $t_0(F_0) = t_0(M)$, so as in the classical case, we can continue building the resolution by replacing $M$ with $\ker(F_0 \to M)$.
\end{proof}

Unlike a minimal resolution, it is not easy to read out the $\Tor$ groups from a degree-minimal flat resolution. However, some inequalities relating the $\Tor$ groups with each component of such a resolution will suffice for our purposes.
\begin{proposition}\label{prop:degreeres}
Let $\bF_{\bullet} \rightarrow M \rightarrow 0$ be a degree-minimal flat resolution of an $A$-module $M$. For all $i \geq 0$, we have
$$ t_0(F_i) \leq \max(t_0(M), t_1(M), \ldots, t_i(M)),$$ 
and 
$$\max(t_0(M),t_1(M), \ldots, t_i(M)) = \max(t_0(F_0), t_0(F_1), \ldots, t_0(F_i)).$$
\end{proposition}
\begin{proof}
   The proof of \cite[Corollary~2.10]{ly17fi} applies.
\end{proof}

\begin{example}\label{exm:koszul2}
    Let $W$ be a finite length $\GL$-representation of positive degree. The Koszul resolution
    \[ \ldots \to \Sym\{W\} \otimes \lw^i\{W\} \to \ldots \to \Sym\{W\} \otimes \lw^2\{W\} \to \Sym\{W\} \otimes W \to \Sym\{W\} \to 0\]
    is a degree-minimal flat resolution of the ground field $k$ over $\Sym\{W\}$. It is rarely a projective resolution as $\lw^i\{W\}$ may not be a projective $\GL$-representation (for example if $W = \bV$). We do not know whether every $A$-module admits a flat cover.
\end{example}

\begin{lemma} \label{lem:torcomparision}
Given $A$-modules $M$ and $N$, and a finite dimensional vector space $V$ over $k$, we have isomorphisms $[\Tor_i^A(M,N)]\{V\} \cong \Tor_i^{A\{V\}}(M\{V\},N\{V\})$ for all $i \geq 0$.
\end{lemma}
\begin{proof}
The result follows from two observations. One, given a $k$-vector space $V$, and $A$-modules $M$ and $N$, the $A\{V\}$-modules $(M \otimes_A N)\{V\}$ and $M\{V\} \otimes_{A\{V\}} N\{V\}$ are naturally isomorphic, and two, the evaluation functor $M \mapsto M\{V\}$ is exact.
\end{proof}

\subsection{\texorpdfstring{$\GL$}{GL}-prime ideals}\label{ss:glprime}
In this section, we prove some preliminary results about $\GL$-prime ideals. Given an $A$-module $M$ and a homogeneous element $m \in M$, we let $\langle m \rangle$ denote the $A$-submodule generated by $m$, i.e., the smallest $\GL$-stable $A$-submodule of $M$ containing the element $m$.  

\begin{definition}
A $\GL$-ideal $\fp$ of a $\GL$-algebra $A$ is \textit{$\GL$-prime} if for all $\GL$-ideals $\fa, \fb$ with $\fa \fb \subset \fp$, either $\fa$ or $\fb$ is contained in $\fp$. A $\GL$-algebra is a \textit{$\GL$-domain} if the ideal $\langle 0 \rangle$ is a $\GL$-prime. 
\end{definition}

We also have an elemental criterion similar to the classical case.
\begin{lemma}\label{lem:elementalcriterion}
    A $\GL$-ideal $\fp$ is $\GL$-prime if and only if for all $f, g \in A$ such that $f (\sigma g) \in \fp$ for all $\sigma \in \GL$ implies either $f \in \fp$ or $g \in \fp$.
\end{lemma}
\begin{proof}
This is standard.
\end{proof}

\begin{definition}
Let $M, N$ be $\GL$-representations and $W$ be a direct summand of $\bV$. A homogeneous element $m \in M$ \textit{belongs to} $M\{W\}$ if the map $M\{\bV\} \to M\{W\} \to M\{\bV\}$ acts by identity on $m$. Two homogeneous elements $m \in M$ and $n \in N$ are \textit{disjoint} if there exists a decomposition $\bV = W_1 \oplus W_2$ such that $m$ belongs to $M\{W_1\}$ and $n$ belongs to $N\{W_2\}$.
\end{definition}

\begin{lemma}\label{lem:disjointnzd}
Let $M$ be an $A$-module. Assume $a \in A$ and $m\in M$ are disjoint elements such that $am=0$. Then $\langle a\rangle \langle m \rangle = 0$.
\end{lemma}
\begin{proof}
Let $\bV = W_1 \oplus W_2$ be a decomposition such that $a$ belongs to $A\{W_1\}$ and $m$ belongs to $M\{W_2\}$. 
For $j =1,2$, let $p_j \colon \bV \to W_j$ and $i_j \colon W_j \to \bV$ be the natural projection and inclusion respectively.

It suffices to show that $a(gm) = 0$ for all $g \in \GL$. Consider the element $e = gi_2p_2 + i_1p_1$ in $\End(\bV)$. We have,
\begin{displaymath}
    0 = e(am) = e(a)e(m) = e(i_1p_1(a))e(i_2p_2(m)) = (i_1p_1(a))(gi_2p_2(m)) = a(gm),
\end{displaymath}
where for the third and last equality, we use that $a = i_1p_1(a)$ and $m = i_2p_2(m)$ and for the penultimate equality, we use that $i_2p_2 i_1p_1 = i_1p_1 i_2p_2 = 0$.
\end{proof}

\begin{corollary}\label{cor:disjointgen}
   Let $a, b \in A$ be disjoint elements. We have $\langle ab \rangle = \langle a \rangle \langle b \rangle$.
\end{corollary}
\begin{proof}
   The containment $\langle ab \rangle \subset \langle a \rangle \langle b \rangle$ clearly holds. The reverse containment follows by passing to $A/\langle ab \rangle $ and using the previous Lemma.
\end{proof}

We note another surprising corollary which shows that a $\GL$-stable ideal is $\GL$-prime if and only if it is $\fS_{\infty}$-prime:
\begin{corollary}\label{cor:glprime=sinfprime}
   Let $\fp$ be a $\GL$-prime ideal in an arbitrary $\GL$-algebra $A$. Then $\fp$ is also $\fS_{\infty}$-prime. 
\end{corollary}
\begin{proof}
    Assume $f, g$ are such that $f \sigma(g) \in \fp$ for all $\sigma \in \fS_{\infty}$. By applying a suitable permutation on $g$, we may further assume that $f$ and $g$ are disjoint so using the previous corollary, we get that $\langle f\rangle \langle g \rangle = \langle fg \rangle \subset \fp$ which implies that either $f$ is in $\fp$ or $g$ is in $\fp$ by $\GL$-primality of $\fp$.
\end{proof}
 
\begin{lemma}\label{lem:xyimpliesxny}
Assume $x$ and $y$ are elements in $A$ such that $xy=0$. Then there exists $m \geq 1$ such that $\langle x^m \rangle \langle y \rangle = 0$.
\end{lemma}
\begin{proof}
For a fixed $N \in \bN$, we first show that for every $a \in \Dist(\GL_N)$, there exists an $m \in \bN$ such that $x^m ay = 0$. When $a=1$, we can take $m = 1$. Assume that for some $a \in \Dist(\GL_N)$ and $s > 0$, we have $x^n e_{i, j}^{(r)}a(y) = 0$ and $x^n\binom{h_i}{r} a(y) = 0$ for all $r < s$. Let $q = p^r$ be the largest power of $p$ that divides $s$. By applying $e_{i, j}^{(q)}$ to the equation $x^n e_{i, j}^{(s-q)}a (y) = 0$, we get 
\begin{displaymath}
0 = e_{i, j}^{(q)} (x^n (e_{i, j}^{(s-q)}a) (y))
    = \sum_{l=0}^{q} e_{i, j}^{(l)}(x^n) (e_{i, j}^{(q - l)}e_{i, j}^{(s-q)}a) (y) 
    = \sum_{l=0}^{(q)} e_{i, j}^{(l)}(x^n) \binom{s-l}{q - l} (e_{i, j}^{(s-l)}a) (y),  
\end{displaymath}
where we use the generalized Leibniz rule to get the second equality, and the multiplication rule in $\Dist(\GL_N)$ to get the third equality. If we multiply the rightmost expression by $x^n$ and use $x^n (e_{i, j}^{(r)}a) (y) = 0$ for $r < s$, only the $l=0$ term survives, so we get $\binom{s}{q} x^{2n} (e_{i, j}^{(s)}a)(y) = 0$. By our choice of $q$, the binomial coefficient $\binom{s}{q} \ne 0$ in $k$ (say, by Lucas's theorem), and so $x^{2n} (e_{i, j}^{(s)}a)(y) = 0$. A similar computation shows that $x^m (\binom{h_i}{s}a)(y) = 0$ for $m \gg 0$. Since the elements $e_{i, j}^{(s)}$ and $\binom{h_i}{s}$ generate $\Dist(\GL_N)$, by induction, we get that for all $a \in \Dist(\GL_N)$, there exists $n$ such that $x^n a y = 0$, as claimed.

Choose an $N \gg 0$ and $\sigma$ a permutation matrix in $\GL_N$ such that $x$ and $\sigma y$ are disjoint. We can find an element $E \in \Dist(\GL_N)$ such that $\sigma y = Ey$ as the $\GL_N$-submodule and the $\Dist(\GL_N)$-submodule generated by the element coincide. By the previous paragraph, we have $x^m Ey = 0$ for some $m \in \bN$, and so by Lemma~\ref{lem:disjointnzd}, we have $\langle x^m \rangle \langle y \rangle = 0$. 
\end{proof}

\begin{remark}
    The previous lemma was proved in characteristic zero by Sam--Snowden \cite[Proposition~8.6.2]{ss12tca} using the universal enveloping algebra.
\end{remark}
\begin{proposition}\label{prop:glprimeisprimary}
    Every $\GL$-prime ideal is primary. The (ordinary) radical of a $\GL$-prime ideal is a $\GL$-stable prime ideal.
\end{proposition}
\begin{proof}
Let $\fp$ be a $\GL$-prime ideal of $A$. By passing to $A/\fp$, we may assume that $\fp = 0$. Assume now that $xy = 0$, then by Lemma~\ref{lem:xyimpliesxny}, we get that $\langle x^n\rangle \langle y \rangle = 0$, but $A/\fp$ is a $\GL$-domain so either $x^n = 0$ or $y = 0$, which implies that $\fp$ is a primary ideal. The claim about the $\rad(\fp)$ now follows since the radical of a primary ideal is prime.
\end{proof}

\subsection{Torsion and generic modules}\label{ss:torgenmod}
Let $M$ be an $A$-module. A nonzero element $m \in M$ is $\textit{torsion}$ if there exists a nonzero $\GL$-ideal $I \subset A$ such that $Im = 0$. An $A$-module $M$ is \textit{torsion} if all nonzero elements of $M$ are torsion. The \textit{torsion submodule} $T(M)$ of $M$ is the submodule generated by all torsion elements. A module $M$ is \textit{torsion-free} if $T(M) = 0$.

\begin{remark}
A $\GL$-equivariant $A$-module which is torsion-free in our sense need not be torsion-free in the classical sense. For example, let $B = \Sym\{\bV\}/\fm^{[p]}$. The algebra $B$ is torsion-free as a $\GL$-equivariant $B$-module as no nonzero ideal annihilates any element in $B$. However, all non-unit elements of $B$ are torsion in the classical sense since they are nilpotent.
\end{remark}
We let $\Mod_A^{\tors}$ be the category of torsion $A$-modules. The subcategory $\Mod_A^{\tors} \subset \Mod_A$ is a localizing subcategory. We let the \textit{generic category} $\Mod_A^{\gen}$ denote the Serre quotient $\Mod_A/\Mod_A^{\tors}$. The canonical \textit{quotient} functor $T \colon \Mod_A \to \Mod_A^{\gen}$ is exact and its right adjoint $\cS$, the \textit{section} functor, is left exact.

\section{Basic Results about \texorpdfstring{$S$}{S}-modules}\label{s:basicresults}
For the remainder of this paper, we let $S = \Sym\{\bV\}$ be the $\GL$-algebra with maximal ideal $\fm$. Recall that an ``$S$-module" means a ``$\GL$-equivariant $S$-module", and ``finitely generated" means ``generated by the $\GL$-orbit of finitely many elements".
\subsection{Noetherianity} \label{ss:noetherianity}
$\GL$-algebras gained prominence since $\GL$-equivariant modules over them are locally noetherian in many examples. In fact, an important question in the theory is whether this property holds for all \textit{finitely generated} $\GL$-algebras. Noetherianity results imply uniformity properties for equations defining natural families of varieties \cite{sno13delta}, and for free resolutions of natural modules \cite{nss16deg2}. For $S$-modules over $k$, the noetherianity property easily follows from Cohen's theorem \cite{co67laws}.

\begin{theorem}\label{thm:locnoeth} 
Let $M$ be a finitely generated $S$-module. Then every $\GL$-stable submodule is also finitely generated.
\end{theorem}
\begin{proof}
Any finitely generated $\GL$-equivariant $S$-module is a quotient of finite direct sums of $P_{\lambda} \coloneqq S \otimes \Div^{\lambda}\{\bV\}$. Furthermore, the module $P_{\lambda}$ is a submodule of $P_n \coloneqq S \otimes \bV^{\otimes n}$, where $n = |\lambda|$ as $\Div^{\lambda}$ is the $\fS_{\lambda}$-invariant subspace of $\bV^{\otimes n}$. So it suffices to show that $P_n$ is $\GL$-noetherian. We identify $\fS_{\infty}$ with the subgroup of permutation matrices in $\GL$. Under this identification, $P_n$ is an $\fS_{\infty}$-equivariant $S$-module and it is finitely generated as $\bV^{\otimes n}$ is a finite length $\fS_{\infty}$-representation. So by Cohen's theorem \cite[\S~2]{co67laws} (see also \cite[Corollary~6.16]{nr19fioi}), the module $P_n$ is $\fS_{\infty}$-noetherian, and so $\GL$-noetherian as well.
\end{proof}

\subsection{Proof of Theorem~\ref{thm:introspec}}\label{ss:glspecofS}
We now concretely analyze the $\GL$-spectrum of $S$ using the results of Section~\ref{ss:glprime}. 

\begin{lemma}\label{lem:x1gen}
Let $n = a_0 + a_1 p + a_2 p^2 + \ldots + a_j p^j$ be the base $p$ expansion of $n$. The $\GL$-ideal generated by $x_1^n$ is $\fm^{a_0} (\fm^{[p]})^{a_1} \ldots (\fm^{[p^j]})^{a_j}$.
\end{lemma}
\begin{proof}
Let $I = \fm^{a_0} (\fm^{[p]})^{a_1} \ldots (\fm^{[p^j]})^{a_j}$. Clearly, $\langle x_1^n \rangle \subset I$ since $x_1^n \in I$. Using the Steinberg tensor product theorem, the $\GL$-representation $I_n$ is the irreducible representation of $\GL$ with highest weight $n$. So $x_1^n$ generates $I_n$, which implies $\langle x_1^n \rangle = I$, as required.
\end{proof}
\begin{lemma}\label{lem:mongen}
Let $n_1, n_2, \ldots, n_l \in \bN$ and write $n_i = a_{i,0} + a_{i,1} p + a_{i,2} p^2 + \ldots + a_{i,j} p^j$ in base $p$. The $\GL$-ideal generated by $x_1^{n_1} x_2^{n_2} \ldots x_l^{n_l}$ is 
$\fm^{b_0} (\fm^{[p]})^{b_1} \ldots (\fm^{[p^j]})^{b_j}$ where $b_m \coloneqq a_{1, m} + a_{2, m} + \ldots + a_{l, m}$ for all $m \geq 0$.
\end{lemma}
\begin{proof}
This follows from Lemma~\ref{lem:x1gen} and Corollary~\ref{cor:disjointgen}.
\end{proof}

\begin{remark}
    Doty \cite{doty89sub} has classified the subrepresentation lattice of $\Sym^d\{\bV\}$ for all $d$. Each submodule gives a distinct $\GL$-stable ideal generated in degree $d$. For example, in characteristic $2$, the ideals generated in degree $10$ are
    \[\fm^{[8]}\fm^{[2]}, (\fm^{[4]})^2\fm^{[2]}, \fm^{[8]}\fm^2, (\fm^{[4]})^2 \fm^2, \fm^{10}.\]
    We caution, however, that $\GL$-stable ideals need not be generated in a single degree, as in $\fm^{10} + \fm^{[8]}$; 
   see \cite{perl24ideals} for a classification.
\end{remark}

We can now compute the $\Spec_{\GL}(S)$ as promised in the introduction.
\begin{proof}[Proof of Theorem~\ref{thm:introspec}]
It is clear that the ideal $(0)$ and $\fm$ are $\GL$-prime since they are prime (in the usual sense). Now fix $q = p^r$. We show that $\fm^{[q]}$ is $\GL$-prime using Lemma~\ref{lem:elementalcriterion}. So assume we have elements $f, g$ such that $f \sigma g \in \fm^{[q]}$ for all $\sigma \in \GL$, and assume $f$ only uses the variables $x_1, x_2, \ldots, x_n$. Take $\sigma$ to be a permutation such that $\sigma g$ only uses the variables $x_i$ with $i > n$. It is now easy to see that either $f \in \fm^{[q]}$ or $\sigma g \in \fm^{[q]}$ (and in the latter case, $g \in \fm^{[q]}$), as required.

We now show the converse that every $\GL$-prime ideal is of the form described in the statement of the theorem. Let $\fp$ be an arbitrary nonzero $\GL$-prime ideal. The only nonzero $\GL$-stable prime ideal of $S$ is the maximal ideal $\fm$. So $\sqrt{\fp} = \fm$ by Proposition~\ref{prop:glprimeisprimary}. Let $n$ be the minimal integer such that $x_1^n \in \fp$. By Lemma~\ref{lem:x1gen}, we see that $\fm^{a_0}(\fm^{[p]})^{a_1} \ldots (\fm^{p^j})^{a_j}$ is contained in $\fp$, where $n = a_0 + a_1 p + \ldots + a_j p^j$ is the base $p$ expansion of $n$. The primality of $\fp$ and minimality of $n$ implies that $n$ must be a power of $p$, say $q$. Therefore $\fp$ contains $\fm^{[q]}$. Assume they are not equal. Let $f \in \fp \setminus \fm^{[q]}$. Since all $\GL$-ideals of $S$ are monomial, we may assume $f$ is a monomial of the form
    $x_1^{\lambda_1}x_2^{\lambda_2}\ldots x_j^{\lambda_j}$ with $\lambda$ a non-increasing sequence and $\lambda_1 < q$. By Lemma~\ref{lem:mongen}, the ideal $\fm^{b_0} (\fm^{[p]})^{b_1} \ldots (\fm^{[p^j]})^{b_j}$ for appropriate integers $b_i$ with $p^j < q$, is contained in $\fp$. By primality of $\fp$, we see that $\fm^{[p^i]}$ is contained in $\fp$, with $p^i < q$, but this contradicts the minimality of $q$. Therefore, $\fp = \fm^{[q]}$, as required.
\end{proof}
For the remainder of this paper, we let $_{\GL}\sqrt{I}$ denote the $\GL$-radical of $I$, i.e., the sum of all $\GL$-ideals $J$ such that $J^N \subset I$ for $N \gg 0$.
\begin{corollary}\label{cor:glrad}
    Let $I$ be a nonzero $\GL$-ideal of $S$. We have $_{\GL}\sqrt{I} = \fm^{[q]}$ for some $q = p^r$.
\end{corollary}
\begin{proof}
We can use standard arguments to show that the $\GL$-radical of $I$ is the intersection of all $\GL$-prime ideals containing $I$. The result now follows from Theorem~\ref{thm:introspec}.
\end{proof}
\begin{remark}\label{rmk:glprimesminors}
    Given an arbitrary $\GL$-algebra $A$, it is natural to wonder whether the Frobenius power of a $\GL$-prime ideal $\fp$ is also $\GL$-prime. We have been able to prove this for ideals of minors in $\Sym\{\bV^{\oplus n}\}$; our proof is not elementary and uses Gr\"obner theory. At the same time, it is not clear how to obtain every $\GL$-prime ideal even for $\Sym\{\bV \oplus \bV\} = k[x_i, y_i \vert i \in \bN]$ in characteristic $2$. For example, the ideal generated by all permutations of $x_1^2, y_1^2, x_1y_1,$ and $x_1y_2+x_2y_1$ is $\GL$-prime.
\end{remark}

\subsection{Torsion modules}\label{ss:torsionSmod}
In characteristic zero, a finitely generated $S$-module is torsion if and only if it has finite length. This fails in positive characteristic: the module $S/\fm^{[p]}$ is evidently a torsion $S$-module, but is supported in infinitely many degrees (the element $x_1 x_2 \ldots x_n$ is in degree $n$), making the torsion category harder to analyze in our setting. 

\begin{lemma}\label{lem:locannihiltor}
Fix $q = p^r$, with $r>0$. Let $M$ be a finitely generated $\Sm$-module. Then $M$ is a torsion $\Sm$-module if and only if every element of $M$ is annihilated by a power of the $\GL$-ideal $\fm^{[q/p]}$. Therefore, the module $M$ is torsion if and only if every element of $M$ is annihilated by the monomial $x^{q/p}_1 x^{q/p}_2 \ldots x^{q/p}_n$ for sufficiently large $n$.
\end{lemma}
\begin{proof}
An element $m \in M$ is torsion if and only if $Im=0$ for some nonzero ideal $I \subset \Sm$. The ideal $I$ will contain a large enough power of $_\GL\sqrt{I}$ since $\Sm$ is $\GL$-noetherian. Therefore $(\fm^{[p^i]})^n \subset I$ for some $p^i < q$ and $n \gg 0$ by Corollary~\ref{cor:glrad}, and in particular $(\fm^{[q/p]})^n \subset I$, proving the first part. The second part now follows by combining the first part with Lemma~\ref{lem:mongen}.
\end{proof}

\begin{lemma}
The Artin--Rees lemma holds for the ideals $\fm^{[q]} \subset S$.
\end{lemma}
\begin{proof}
Let $I = \fm^{[q]}$. The Rees algebra of $I$ is a quotient of $S' = \Sym\{\bV\} \otimes \Sym\{\bV^{(r)}\}$. The action of $\fS_{\infty}$ on $\bV^{(r)}$ is the same as the action of $\fS_{\infty}$ on $\bV$, so $S'$ is $\GL$-noetherian by Cohen's theorem \cite[\S~2]{co67laws} (see also \cite[Corollary~6.16]{nr19fioi}) and therefore, so is the Rees algebra. The proof from the classical case now implies the lemma.
\end{proof}
\begin{remark}
It is remarkable that we are able to prove the Artin--Rees lemma for $\GL$-prime ideals in $S$ using known noetherianity results. We currently do not know how to prove the above result for the ideal of $2 \times 2$-minors in $S \otimes S$ over fields of positive characteristic; in characteristic zero, the Rees algebra associated with any ideal in $S \otimes S$ is known to be noetherian since they are \textit{bounded} (see \cite[Corollary~4.20]{ss19gl2}). 
\end{remark}

We can now show that property (Inj) from Section~\ref{ss:torgenmod} holds for certain subcategories of interest.
\begin{proposition}\label{prop:propinj}
The following subcategories satisfy property (Inj), i.e., injectives in the subcategory remain injective in the ambient category:
\begin{enumerate}
    \item $\Mod_{\Sm}^{\tors} \subset \Mod_{\Sm}$, 
    \item $\Mod_S[(\fm^{[q]})^{\infty}] \subset \Mod_S$, and
    \item $\Mod_S^{\tors} \subset \Mod_S$.
\end{enumerate}
\end{proposition}
\begin{proof}
By Lemma~\ref{lem:locannihiltor}, the category $\Mod_{\Sm}^{\tors}$ is the full subcategory of $\Mod_{\Sm}$ on the objects locally annihilated by $\fm^{[q/p]}$. Since the Artin--Rees lemma holds for these ideals, injective objects in $\Mod_{\Sm}^{\tors}$ remain injective in $\Mod_{\Sm}$ by \cite[Corollary~4.19]{ss19gl2}. 

The argument in the previous paragraph shows that $\Mod_S[(\fm^{[q]})^{\infty}] \subset \Mod_S$ also satisfies property (Inj). 

For the last part, we first show that an arbitrary torsion $S$-module can be embedded into an injective in $\Mod_S$ that is torsion. Let $U$ be such a module. We can write $U$ as the direct limit of its finitely generated submodules, say $\{U_{\lambda}\}$. Being finitely generated, each $U_\lambda$ lies in the subcategory $\Mod_S[(\fm^{[q]})^{\infty}]$ for sufficiently large $q$, so the injective envelope $I_{\lambda}$ of $U_{\lambda}$ in $\Mod_S[(\fm^{[q]})^{\infty}]$ will also be injective in $\Mod_S$ by the previous paragraph. The direct limit $I$ of the various $I_{\lambda}$ will be a torsion $S$-module that contains $U$; it is also injective in $\Mod_S$ since the direct limit of injective modules in a locally noetherian Grothendieck abelian category is injective. 

Now, given an arbitrary injective $I'$ in $\Mod_S^{\tors}$, we can embed $I'$ into a torsion $I$ which is injective in $\Mod_S$ by the previous paragraph. This map splits in $\Mod_S^{\tors}$ as $I'$ is injective, so $I'$ is a direct summand of $I$ and hence injective in $\Mod_S$ as well.
So, $\Mod_S^{\tors} \subset \Mod_S$ also satisfies property (Inj).
\end{proof}

\section{Embedding Theorem}\label{s:snowdenshift}
The goal of this section is to study the generic category of $\Mod_S^{\gen}$, by extending a recent result of Snowden \cite{sno21stable} about torsion-free modules over polynomial $\GL$-algebras to all characteristics.

\subsection{Tensor-disjoint weights} \label{ss:flatwgts} 
In this section, we prove some technical results about weight vectors in tensor products. 

\subsubsection{Flat weights}
Snowden \cite[Section~3]{sno21stable} proved Theorem~\ref{thm:embeddingintro} in characteristic zero by exploiting the presence of weight vectors of weight $(1, 1, \ldots, 1)$ in any polynomial representation of $\GL$. This phenomenon is unique to the infinite general linear group (for example, the $\GL_n$-representation $\Sym^d\{k^n\}$ has such a vector if and only if $n \geq d$) and occurs only in characteristic zero (the Frobenius power of any $\GL$-representation will not have such a weight). To overcome this difficulty, we introduce the notion of a \textit{flat} weight.
\begin{definition}\label{def:minflatwgt}
A weight $\lambda$ is \textit{flat} if each component of $\lambda$ is a power of $p = \chark(k)$ or $0$. The \textit{$p$-magnitude} of the weight $\lambda$, denoted $\pmag(\lambda)$, is the sequence $(n_0, n_1, \ldots, n_r, \ldots )$ where $n_i = \#\{\lambda_j | \lambda_j = p^i\}$. 
\end{definition}
We define a total ordering $<$ on the set of flat weights using $\pmag$: given flat weights $\lambda$ and $\mu$, we have $\lambda < \mu$ if and only if $\pmag(\lambda) > \pmag(\mu)$ in the lexicographic ordering on sequences of integers. We transport this definition to weight vectors in the obvious way.
\begin{definition}\label{def:minflatvec}
Given a $\GL$-representation $W$, a nonzero weight vector $w$ is \textit{flat} if $\wgt(w)$ is flat.
The weight vector $w \in W$ is \textit{minimally flat} if $\wgt(w)$ is minimal among all flat weights occurring in $W$ under the total order $<$ defined on flat weights.
\end{definition}

\begin{example}
Let $W = \Sym^p$. The elements $x_1^p$ and $x_1 x_2 x_3 \ldots x_p$ are both flat; $x_1x_2\ldots x_p$ is minimally flat in $\Sym^p$.
\end{example}

\begin{proposition}\label{prop:minflatirrep}
 Let $\mu$ be a partition with $\mu = \mu^0 + p \mu^1 + \ldots + p^r \mu^r$ being the $p$-adic expansion of $\mu$ into $p$-restricted partitions. 
 Assume $v$ is a minimally flat weight vector in $L_{\mu}$. Then $\pmag(\wgt(v)) = (|\mu_0|, |\mu_1|, \ldots, |\mu_r|)$. 
\end{proposition}
\begin{proof}
    Given a $p$-restricted partition $\lambda$, the $\GL$-representation $L_{\lambda}$ has a weight vector of weight $1^{|\lambda|}$ \cite[Theorem~1]{pre87wts}, which is clearly minimally flat in $L_{\lambda}$. The result follows by the Steinberg tensor product theorem.
\end{proof}

\subsubsection{The shift functor}
We let $G(n)$ denote the subgroup of $\GL$ consisting of block matrices of the form
$$\begin{pmatrix}
    \id_n & 0 \\
    0 & * 
\end{pmatrix}$$
where the top block is the identity matrix of size $n \times n$. 
The subgroup $G(n)$ is isomorphic to $\GL$ with isomorphism given by $A \mapsto \begin{pmatrix}
    \id_n & 0 \\ 0 & A
\end{pmatrix}$. 

\begin{definition}
The \textit{shift} $\nSh_n(W)$ of a $\GL$-representation $W$ is the $\GL$-representation with the same underlying space as $W$ and $\GL$ acting via the self-embedding $\GL \cong G(n) \subset \GL$ given above.
\end{definition}

The next lemma shows that the shift functor does not change a representation in its top degree.
\begin{lemma}[Lemma~14 in \cite{dra19top}]\label{lem:shiftlowers}
Assume $V$ is a homogeneous polynomial representation of degree $d$. The $\GL$-representation $\nSh_n(V) \cong V \oplus W$ where $W$ is a polynomial representation of degree $< d$.
\end{lemma}

\begin{remark}
Given a $\GL$-algebra $A$ and an $A$-module $M$, the $\GL$-representation $\nSh_n(A)$ has the canonical structure of a $\GL$-algebra and $\nSh_n(M)$ is canonically a $\nSh_n(A)$-module as we are merely restricting to a subgroup of $\GL$. 
\end{remark}

We will only use $\nSh_n$ in this section. In particular, what we refer to as the shift functor in Section~\ref{s:nagpalshift} is not the same as this one (but they are closely related).

\subsubsection{Technical results about flat weights}
In what follows, we let the $\supp(\lambda)$ be the support of a weight $\lambda$, which is the set of natural numbers $i$ such that $\lambda_i \ne 0$. 

We can now prove the technical result about tensor-disjoint weight vectors in any subrepresentations alluded to in Section~\ref{ss:embedintro}.
\begin{proposition}\label{prop:disjointweights}
Let $V$ and $W$ be $\GL$-representations of positive degrees $n$ and $m$ respectively. Let $U \subset V \otimes W$ be a nonzero subrepresentation. Let $u$ be a minimally flat weight vector in $U$ with $u = \sum_{i=1}^s v_i \otimes w_i$ where the $\{v_i\}$ and $\{w_i\}$ are linearly independent weight vectors of $V$ and $W$ respectively. For all $i \geq 1$, the weights $\wgt(v_i)$ and $\wgt(w_i)$ have disjoint support.
\end{proposition}
\begin{proof}
    Assume the weight of $u$ is $\lambda$ with $\lambda_i = 0$ for $i > r+1$. Let $v_i$ and $w_i$ have weight $\mu^i$ and $\nu^i$ respectively. We can swap any component of the weight of $\lambda$ with its $(r+1)$-st component by applying a suitable permutation, so the lemma follows if we can show that for all $i$, either $\mu^i_{r+1} = 0$ or $\nu^i_{r+1} = 0$. 
    
    Assume $\lambda_{r+1} = p^l$.
    Let $U' \subset U$ be the $G(r)$-representation generated by $u$. The vector $u$ has weight $(p^l)$ for the action of $G(r)$, so $U'$ is a degree $p^l$ polynomial $G(r)$-representation. 
    
    We first show that the $G(r)$-representation $U'$ is isomorphic to the $l$-fold Frobenius twist of the standard representation $k\langle e_{r+1}, e_{r+2}, \ldots, e_N, \ldots \rangle$ of $G(r)$. Assume for the sake of contradiction that it is not. Then $U'$ contains a flat $G(r)$-weight vector $u'$ with all the components of its weights strictly smaller than $p^l$ by Proposition~\ref{prop:minflatirrep}. Furthermore, since the action of top $r$-dimensional torus in the ambient $\GL$ acts on $u'$ and $u$ are the same, $u'$ is also a flat weight vector for the ambient $\GL$-action.
    However, the weight of $u'$ is strictly smaller than $\wgt(u)$ contradicting the choice of $u$, so $U'$ is $l$-fold Frobenius twisted, as claimed.

    Put $\nSh_n(V) \cong \tilde{V}_0 \oplus \tilde{V}_1 \oplus \ldots \oplus \tilde{V}_n$ where $\tilde{V}_i$ are degree $i$ polynomial representations of $\GL$ and similarly, $\nSh_n(W) \cong \tilde{W}_0 \oplus \tilde{W}_2 \oplus \ldots \oplus \tilde{W}_m$. Then $\nSh_n(V\otimes W) \cong \bigoplus_{0 \leq i \leq n, 0 \leq j \leq m} ({\tilde{V}_i \otimes \tilde{W}_j})$, and set $Z \subset \nSh_n(V \otimes W)$ to be the direct sum of $\tilde{V_i} \otimes \tilde{W_j}$ with both $i$ and $j$ strictly positive. We have a canonical surjection of $\GL$-representations $\nSh_n(V \otimes W) \to Z$. Under the composite map $U' \to (V \otimes W)\vert_{G(r)} \to Z$, the image of $u$ is the sum of all $v_j \otimes w_j$ with the $G(r)$-weight of both $v_j$ and $w_j$ nonzero.
    
   Assume now that there exists an $i$ such that both $\mu^i_{r+1} \ne 0$ and $\nu^i_{r+1} \ne 0$. Then the image of $u$ under the composite map is nonzero by the previous paragraph. However, this contradicts the fact that $\Hom_{\GL}(\bV^{(l)}, T \otimes T') = 0$ for all positive degree $\GL$-representations $T, T'$ by \cite[Theorem~2.12]{fs97coh}.
\end{proof}

We now prove a variant of Snowden's key lemma \cite[Lemma~3.2]{sno21stable}.

\begin{lemma}
    Let $V$ and $W$ be homogeneous polynomial representations of degree $n$ and $m$ with $V$ irreducible. Let $U \subset V \otimes W$. Then there exists a weight vector $u = \sum_{i=0}^r v_i \otimes w_i$ such that 
\begin{itemize}
    \item for all $i \in [r]$, both $v_i$ and $w_i$ are weight vectors with $\wgt(v_i) + \wgt(w_i) = \wgt(u)$; 
   \item for all $i \in [r]$, the weights $\wgt(v_i)$ and $\wgt(w_i)$ have disjoint support; 
    \item $\supp(\wgt(u)) = [a + b]$ for some $a, b > 0$;
    \item $\supp(\wgt(v_1)) = [a]$; and
    \item $\supp(\wgt(v_i)) \ne [a]$ for $i > 1$.
\end{itemize}
\end{lemma}
\begin{proof} 
    By the previous proposition, we can find $u = \sum_{i=1}^r v_i \otimes w_i$ satisfying the first two conditions. By collecting terms if necessary, we can further assume that $\{v_i\}_{i\in [r]}$ is a linearly independent set. Finally, by re-indexing and applying a permutation in $\GL$, we can also assume that
    \begin{itemize}
        \item $\supp(\wgt(u)) = [a+b]$,
        \item $\supp(\wgt(w_1)) = [a]$ and $\supp(\wgt(v_1)) = [a+b]\setminus [a]$,
        \item $\wgt(v_i) = \wgt(v_1)$ for $i \leq l$, and
        \item $\wgt(v_i) \ne \wgt(v_1)$ for $i > l$.
    \end{itemize}
Note that since $\wgt(v_i)$ and $\wgt(w_i)$ have disjoint supports, the last two conditions above are equivalent to $\supp(\wgt(v_i)) = [a+b] \setminus [a]$ for $i \leq l$ and $\supp(\wgt(v_i)) \ne [a+b] \setminus [a]$ for $i > l$.
    
Let $V' \subset V$ be the $G(a)$-representation generated by $v_1, v_2, \ldots, v_l$. 
Identifying $G(a)$ with $\GL$, we have an isomorphism $V' \cong V$ by Lemma~\ref{lem:shiftlowers} so $V'$ is irreducible. Furthermore, span$\langle w_1 \rangle \subset W$ is the trivial representation of $G(a)$, and by assumption $V' \ncong k \langle w_1 \rangle$. Therefore, using the Jacobson Density Theorem \cite[Theorem~3.2.2(ii)]{eetal11rep}, we can find an element $a \in k[G(a)]$ such that  $a(v_1) = v_1$, $aw_1 = w_1$, and $a(v_i) = 0$ for $2 \leq i \leq l$. Write $\nSh_a(V) \cong V' \oplus Z$. We have $
\nSh_a(V \otimes W) \cong (V' \otimes \nSh_a(W)) \oplus (Z \otimes \nSh_a (W))$. Note that for $i > l$, the elements $v_i \otimes w_i$ lie in the second summand, therefore $a(v_i \otimes w_i) \in Z \otimes \nSh_n(W)$ for $i > l$. So $au = v_1 \otimes w_1 + a(\sum_{i= l+1}^r v_i \otimes w_i) \ne 0$. For a suitable permutation matrix $\sigma \in \GL$, one easily checks that the element $\sigma au$ satisfies all the conditions, as required.
\end{proof}
\begin{remark}\label{rmk:bkprest}
The results of this subsection are somewhat delicate, as illustrated in \cite[Remark~3.5]{bk99tensor}. Brundan--Kleshchev observe that the socle of the tensor product of two irreducible $p$-restricted representations need not contain any $p$-restricted irreducibles. Consequently, one cannot rely on characteristic-zero arguments, such as the existence of weight vectors of the form $1^n$, to obtain tensor-disjoint weight vectors.
\end{remark}

\subsection{Proof of Theorem~\ref{thm:embeddingintro}}\label{ss:embedding}

\begin{definition}
    A $\GL$-algebra $A$ over $k$ is \textit{integral} if the underlying $k$-algebra $A$ is an integral domain in the usual sense.
    A $\GL$-algebra $A$ is \textit{polynomial} if $A = A_0 \otimes \Sym\{W\}$ for a polynomial representation $W$, where $A_0$ is the degree zero component of $A$.
\end{definition}

For this section, we recall a well-partial order $\prec$ on finite length $\GL$-representations (see, for example, \cite[\S2.2]{dra19top}). Given polynomial representations $V$ and $W$, we say $V \prec W$ if $V$ and $W$ are non-isomorphic and in the largest degree where $V_e \not\cong W_e$, there is a surjection $W_e \to V_e$. 

The key lemma about flat weights from the previous section is used in the proof of the next proposition. This proposition is essentially a ``linearization" of a technical result about $\GL$-varieties proved by Draisma \cite{dra19top} (see also \cite[Theorem~4.1]{bdes21geo}) which states that given a closed subset $Z$ of an affine $\GL$-variety $\Spec(\Sym\{W\})$, a Zariski open subset in $\nSh_n(Z)$ embeds into $\Spec(\Sym\{W'\})$ where $W' \prec W$.
\begin{proposition} \label{prop:reduction}
Let $A$ be an integral $\GL$-algebra, let $\lambda$ be a partition, let $F$ and $M$ be $A$-modules, and let $P$ be an extension
\begin{displaymath}
0 \to F \to P \to A \otimes L_{\lambda}  \to 0.
\end{displaymath}
Assume further that we have a surjection $P \to M$.
Then at least one of the following holds:
\begin{enumerate}
\item $\ker(P \to M) \subset F$
\item There exists $n \ge 0$ and a non-zero $\GL$-invariant element $f \in \nSh_n(A)$ such that the natural map $\nSh_n(P) \to \nSh_n(M)$ restricts to a surjection $Q[1/h] \to \nSh_n(M)[1/h]$, where $Q \subset \nSh_n(P)$ is the preimage of the first summand of 
$$\nSh_n(A \otimes L_{\lambda}) = (\nSh_n(A) \otimes W) \oplus (\nSh_n(A) \otimes L_{\lambda}),$$
\end{enumerate}
where $\nSh_n(L_{\lambda}) \cong W \oplus L_{\lambda}$ as in Lemma~\ref{lem:shiftlowers}.
\end{proposition}

\begin{proof}
Proof of \cite[Proposition~3.1]{sno21stable} applies almost verbatim.
\end{proof}

This next result states that a module is flat after shifting and localizing; it is obtained by iterating the conclusion of the previous proposition. 
\begin{proposition} \label{prop:genshift}
Let $A$ be an integral $\GL$-algebra and let $M$ be a finitely generated $A$-module. There exist $n \geq 0$ and a nonzero $\GL$-invariant function $f \in \nSh_n(A)$ such that $\nSh_n(M)[1/f]$ is a semi-induced $\nSh_n(A)[1/f]$-module.
\end{proposition}
\begin{proof}
We adjust the proof of the same result of Snowden in characteristic zero \cite[Theorem~3.3]{sno21stable}. Say an $A$-module is \textit{good} if the conclusion of the theorem holds for it. For a finite length $\GL$-representation $W$, an $A$-module $M$ is an \textit{essential quotient} of $W$ if $M$ is a quotient of a semi-induced $A$-module $P$ with $\Tor_0^A(P, k) = W$.
Consider the following statement, for a $\GL$-representation $W$:

$\Psi(W)$ -- If $A$ is an integral $\GL$-algebra and $M$ is an essential quotient of $W$, then $M$ is good.

We prove $\Psi(W)$ for all finite length $\GL$-representations by induction on the partial order $\prec$. The base case, when $V=0$ is trivial. Assume we know $\Psi(V)$ for all $V \prec W$. We now prove it for $W$.

Let $A$ be an integral $\GL$-algebra and let $M$ be a quotient of $P$ with $\Tor_0^A(P, k) = W$. Let $L_{\lambda}$ be an irreducible quotient of $W$ of maximal degree $d$. We obtain a surjection $P \to A \otimes L_{\lambda}$ by composing the surjection $P \to P/P^{<d} \cong A \otimes W_d$, with the surjection $A \otimes W_d \to A \otimes L_{\lambda}$.

We apply Proposition~\ref{prop:reduction} to the short exact sequence
\begin{displaymath}
    0 \to F \to P \to A \otimes L_{\lambda} \to 0
\end{displaymath}
and the surjection $P \to M$; the $A$-module $F$ is also semi-induced by Corollary~\ref{cor:semises}.

Suppose case (1) holds. Then $M$ is an extension of $F/N$ for some submodule $N \subset F$, and $A \otimes L_{\lambda}$. The $A$-module $F$ is an essential quotient of $W/L_{\lambda}$, and $W/L_{\lambda} \prec W$, so $\Psi(W/L_{\lambda})$ holds which means $F/N$ is good, and therefore so is $M$.

Suppose instead case (2) holds. Then there exists $n \ge 0$ and a $\GL$-invariant function $f \in \nSh_n(A)$ such that the natural map $\nSh_n(P) \to \nSh_n(M)$ restricts to a surjection $Q[1/f] \to \nSh_n(M)[1/f]$, where $Q \subset \nSh_n(P)$ is the extension of $\nSh_n(F)$ and $\nSh_n(A) \otimes V$ with $V$ being the 
degree $<d$ subrepresentation of $\nSh_n(L_{\lambda})$ by Lemma~\ref{lem:shiftlowers}.
Set $U \coloneqq \Tor_0^{\nSh_n(A)[1/f]}(Q[1/f], k)$. Then $U$ is an extension of $\nSh_n(W/L_{\lambda})$ and $V$, so by Lemma~\ref{lem:shiftlowers} we see $U \prec W$. By induction, the statement $\Psi(U)$ holds and so the $\nSh_n(A)[1/f]$-module $\nSh_n(M)[1/f]$ is good being the quotient of the semi-induced $\nSh_n(A)[1/f]$-module $Q[1/f]$ which readily implies that $M$ is also a good $A$-module, as required.
\end{proof}

\begin{proof}[Proof of Theorem~\ref{thm:embeddingintro}]
Let $W$ be a finite length polynomial representation of $\GL$ and $A = \Sym\{W\}$. Writing $\nSh_n(W) \cong W \oplus \tilde{W} \oplus U$, where $\tilde{W}$ is purely of positive degree, and $U$ is of degree zero, we have $\nSh_n(A) \cong A\otimes \Sym\{\tilde{W}\} \otimes \Sym\{U\}$. By restricting scalars, the $A$-module $\nSh_n(A)$ is an infinite direct sum of finitely generated induced $A$-modules, so any finitely generated $A$-submodule of $\nSh_n(A)$ lies in a finitely generated induced $A$-module, and in turn, any finitely generated $A$-submodule of a semi-induced $\nSh_n(A)$-module $F$ also lies in a finitely generated semi-induced submodule of $F$.

Given a finitely generated torsion-free $A$-module $M$, by Proposition~\ref{prop:genshift}, there exists $n > 0$ and $\GL$-invariant function $f \in \nSh_n(A)$ such that $\nSh_n(M)[1/f]$ is a semi-induced $\nSh_n(A)[1/f]$-module. Let $F$ be a semi-induced $\nSh_n(A)$-module such that the $\nSh_n(A)[1/f]$-module $\nSh_n(M)[1/f]$ is isomorphic to $F[1/f]$. Since $M$ is torsion-free, the canonical map $M \to \nSh_n(M) \to \nSh_n(M)[1/f]$ is injective, and by scaling this with a power of $f$, we may assume that the image of $M$ lands in $F$. By the previous paragraph, the image of $M$ lands inside a finitely generated semi-induced $A$-submodule $G \subset F$. The map $M \to G$ satisfies the required property.
\end{proof}

\subsection{Torsion-free injectives}\label{ss:Sinj}
In this section, we use the embedding theorem to show that $S \otimes I$ is an injective $S$-module when $I$ is an injective $\GL$-representation.
\begin{lemma}\label{lem:finflatdim}
    Let $M$ be a finitely generated $S$-module, and assume $\Tor_i^S(M, k) = 0$ for some $i > 0$. Then $M$ is semi-induced.
\end{lemma}
\begin{proof}
    The $S$-module $M$ is semi-induced if and only if $\Tor_1^S(M,k)=0$ by
    Proposition~\ref{prop:flatequalssemi}. By an easy induction, it suffices to show that if $\Tor_2^S(M, k) = 0$, then $\Tor_1^S(M, k) = 0$. Assume not. Then using Lemma~\ref{lem:torcomparision} for $n \gg 0$ we have $\Tor_1^S(M\{k^n\}, k) \ne 0$ and $\Tor_2^S(M\{k^n\}, k) = 0$. So $M\{k^n\}$ has depth $n-1$ as a $k[x_1, x_2, \ldots, x_n]$-module by the Auslander--Buchsbaum formula. Similarly, the $\GL_{n+1}$-equivariant $k[x_1, x_2, \ldots, x_{n+1}]$-module $M\{k^{n+1}\}$ has depth $n$. By the reasoning in \cite[Lemma~7.3.2]{ss16gl}, the sequence $x_1, x_2, \ldots, x_n$ forms a regular sequence in $M\{k^{n+1}\}$ (recall we are assuming $k$ is infinite). The $k[x_1, x_2, \ldots, x_n]$-module $M\{k^n\}$ is a direct summand of $M\{k^{n+1}\}$, and so the sequence $x_1, x_2, \ldots, x_n$ is also a regular sequence on $M\{k^n\}$, contradicting the fact that the depth of $M\{k^n\}$ is $n-1$. So $\Tor_1^S(M, k)=0$, as required. 
\end{proof}

\begin{corollary}\label{cor:noinjectivemaps}
   Let $F$ be a semi-induced $S$-module generated in degree $< n$ and $V$ be a homogeneous polynomial representation of degree $n$. There are no injective maps $S \otimes V \to F$.
\end{corollary}
\begin{proof}
    Assume $f$ is such an injective map $S \otimes V \to F$. Then the cokernel of $f$ satisfies
    \[\Tor_1^S(\coker(f), k) \cong V \text{ and } \Tor_i^S(\coker(f),k) = 0 \text{ for }i > 1,\]
    contradicting the previous lemma.
\end{proof}

\begin{corollary}\label{cor:semiindsplit}
    Let $F$ be a finitely generated semi-induced $S$-module and let $I$ be an injective polynomial $\GL$-representation. Any injective map $S \otimes I \to F$ splits.
\end{corollary}
\begin{proof}
Recall that the tensor product of two injective polynomial representations is also injective in $\Rep^{\pol}(\GL)$ \cite[Corollary~2.11]{fs97coh}, so $S \otimes I$ is injective as a polynomial representation. For an induced module $S \otimes V$, we have $\ext^i_{S}(S \otimes V, S\otimes I) = \ext^i_{\GL}(V, S\otimes I) = 0$ for $i > 0$. By induction, we also see that $\ext^i(M, S\otimes I) = 0$ for any semi-induced module $M$. 

Given an injection $f \colon S \otimes I \to F$, the cokernel is semi-induced by Lemma~\ref{lem:finflatdim}, so the injection splits as $\ext^1(\coker(f), S \otimes I) = 0$  by the previous paragraph.
\end{proof}
\begin{proposition}\label{prop:injective}
    For an injective $\GL$-representation $I$, the $S$-module $S \otimes I$ is injective. 
\end{proposition}
\begin{proof}
    Let $N = S \otimes I$. 
    By Baer's criterion \cite[Proposition~A.4]{gs18incmon}, it suffices to show that $\ext^1(M, N) = 0$ for \textit{finitely generated} $M$. For a torsion module, we know $\ext^1(T, N) = 0$ by Proposition~\ref{prop:propinj}, so we may further assume that $M$ is \textit{torsion-free}. So assume we have an extension
    \[ 0 \to N \to K \to M \to 0\]
    in $\ext^1(M, N)$. The $S$-module $K$ is also torsion-free, being an extension of two torsion-free modules. By Theorem~\ref{thm:embeddingintro} we have an injection $K \to F$ with $F$ a semi-induced $S$-module. Composing this with the injection $N \to K$, we obtain an injection $N \to F$, which splits by Corollary~\ref{cor:semiindsplit} so in particular the injection $N \to K$ also splits, which implies $\ext^1(M, N) = 0$, as required.
\end{proof}
\begin{remark}\label{rmk:touzecomment}
    Let $A = \Sym\{W\}$ with $W$ a finite length $\GL$-representation. In characteristic zero, Snowden \cite[Corollary~4.11]{sno21stable} showed that $A \otimes \bS_{\lambda}\{\bV\}$ is injective in $\Mod_A$ by first proving that they are injective in $\Mod_A^{\gen}$ using the embedding theorem, and then using the easy fact that the section functor $\Mod_A^{\gen} \to \Mod_A$ preserves injective objects. The analogous result fails for arbitrary $W$ in characteristic $p$.
    
    Let $W = \Sym^2$. In this case, $A$ is not self-injective, as we now explain. It is easy to see that given an injective object $M$ in $\Mod_A$, the graded pieces of $M$ are injective in $\Rep^{\pol}(\GL)$. But $\Sym^p(\Sym^2)$ is not injective in $\Rep^{\pol}(\GL)$ (see \cite[Section~6]{tou13ringel}).
\end{remark}

\subsection{The generic category}\label{ss:genericSmod}
Snowden \cite{sno21stable} uses the embedding theorem to study the generic categories of $\GL$-algebras. His methods do not directly extend to positive characteristic owing to $\Rep^{\pol}(\GL)$ not being semisimple. However, the algebra $S$ is fairly well behaved --- for example, its graded pieces are injective in $\Rep^{\pol}(\GL)$ --- allowing us to study the generic category of $S$ by making minor modifications to Snowden's techniques.
Throughout, we let $T$ be the exact functor $\Mod_S \to \Mod_S^{\gen}$ and set $J_{\lambda} = T(S \otimes I_{\lambda})$ and $K_{\lambda} = T(S \otimes L_{\lambda})$. We first note a refinement of Theorem~\ref{thm:embeddingintro} when $A = S$. 
\begin{proposition}\label{prop:finalemb}
   Every finitely generated torsion-free $S$-module injects into a finite direct sum of $S \otimes I_{\lambda}$ (with $\lambda$ allowed to vary). Consequently, every finitely generated object in $\Mod_S^{\gen}$ injects into a finite direct sum of $J_{\lambda}$.
\end{proposition}
\begin{proof}
By Theorem~\ref{thm:embeddingintro}, it suffices to show that every finitely generated semi-induced module injects into a finite direct sum of $S \otimes I_{\lambda}$, which is clear since $S \otimes L_{\lambda}$ injects into $S \otimes I_{\lambda}$ and $S \otimes I_{\lambda}$ is injective.
\end{proof}

\begin{lemma}
    The object $J_{\lambda}$ is an indecomposable injective in $\Mod_S^{\gen}$. 
\end{lemma}
\begin{proof}
The $S$-module $S \otimes I_{\lambda}$ is an indecomposable injective, so its image is an indecomposable injective in $\Mod_S^{\gen}$ by \cite[Proposition~4.4]{ss16gl}.
\end{proof}

We note a useful result.
\begin{proposition}\label{prop:cleanupprop}
    Let $\cA$ be a Grothendieck abelian category and assume $\cB$ is a localizing subcategory satisfying property (Inj) with $\Gamma \colon \cA \to \cB$ right adjoint to the inclusion functor, and $\cS \colon \cA/\cB \to \cA$ right adjoint to the exact functor $T \colon \cA \to \cA/\cB$. For any object $M$ in $\cA$, we have a natural four-term short exact sequence:
    \[ 0 \to \Gamma(M) \to M \to \cS(T(M)) \to \rR^1\Gamma(M) \to 0,\] 
    where the map $M \to \cS(T(M))$ is the natural unit map. Furthermore, for all $i \geq 1$, we have isomorphisms $\rR^i\cS(T(M)) \cong \rR^{i+1}\Gamma(M)$.
\end{proposition}
\begin{proof}
This follows from (the proof of) \cite[Proposition~4.6]{ss19gl2}. Note, however, that the proof only shows that $\rR^{i+1}\Gamma(M) \cong \rR^i(\cS \circ T)(M)$, but using the exactness of $T$, we in fact have an isomorphism of functors $\rR^i(\cS \circ T) \cong (\rR^i\cS) \circ T$ whence the claimed result follows.
\end{proof}
\begin{corollary}\label{cor:injsat}
The natural unit map $S \otimes I_{\lambda} \to \cS \circ T (S \otimes I_{\lambda})$ is an isomorphism.
\end{corollary}
\begin{proof}
    The $S$-module $S \otimes I_{\lambda}$ is torsion-free, so $\Gamma(S \otimes I_{\lambda}) = 0$, and $\rR^1\Gamma(S \otimes I_{\lambda}) = 0$ as $S \otimes I_{\lambda}$ is injective by Proposition~\ref{prop:injective}. The result now follows by applying Proposition~\ref{prop:cleanupprop}.
\end{proof}

\begin{lemma}\label{lem:adjunctions}
    For a $\GL$-representation $V$, we have 
    \[\Hom_{\Mod_S^{\gen}}(T(S \otimes V), J_{\lambda}) \cong \Hom_{S}(S \otimes V, S \otimes I_{\lambda}) \cong  \Hom_{\GL}(V, S \otimes I_{\lambda}).\]
\end{lemma}
\begin{proof}
The first isomorphism follows by using the definition of adjunction combined with the previous corollary. The second isomorphism follows because the forgetful functor $\Mod_S \to \Rep^{\pol}(\GL)$ is left adjoint to the functor $V \mapsto S \otimes V$.
\end{proof}

\begin{corollary}\label{cor:genericmap}
    Let $f \colon K_{\lambda} \to J_{\mu}$.
    \begin{itemize}
        \item If $|\mu| \geq |\lambda|$ and $\mu \neq \lambda$, then $f$ is zero, 
        \item if $|\mu| = |\lambda|$, then $f$ is either injective or zero, and
        \item if $|\mu| < |\lambda|$, then $f$ is not injective.
    \end{itemize}
\end{corollary}
\begin{proof}
   By the previous lemma, we have 
   \[\Hom_{\Mod_S^{\gen}}(K_{\lambda}, J_{\mu}) \cong \Hom_{\GL}(L_{\lambda}, S \otimes I_{\lambda}) = \Hom_{\GL}(S_{|\lambda|-|\mu|} \otimes I_{\lambda})\]
   which is trivially zero when $|\mu| > |\lambda|$. When $|\mu| = |\lambda|$ and $\mu \ne \lambda$, then $\Hom_{\GL}(L_{\lambda}, I_{\mu}) = 0$ since there are no nonzero morphisms from a simple object to the injective envelope of a non-isomorphic simple object in any category of locally finite length. When $\mu = \lambda$, either $f$ maps $L_{\lambda}$ isomorphically onto the socle of $I_{\mu}$ in which case $f$ is injective, or it kills $L_{\lambda}$ in which case $f = 0$. The final case follows from Corollary~\ref{cor:noinjectivemaps}.

\end{proof}
\begin{lemma}\label{lem:kernelissimple}
    Let $H_{\lambda}$ be the intersection of the kernels of all maps from $K_{\lambda}$ to a finite direct sum of $J_{\mu}$ with $|\mu| < |\lambda|$. The object $H_{\lambda}$ is simple.
\end{lemma}
\begin{proof}
    There are no injective maps from $K_{\lambda}$ to a finite direct sum of $I_{\mu}$ with $|\mu| < |\lambda|$ by Lemma~\ref{lem:adjunctions} and Corollary~\ref{cor:noinjectivemaps}, so the module $H_{\lambda}$ is nonzero.
    Assume $H_{\lambda}$ has a nonzero subobject $N$. Then we can find an injection $i \colon K_{\lambda}/N \to \bigoplus J_{\mu}$ by Proposition~\ref{prop:finalemb}. Composing $i$ with the canonical surjection $\pi \colon K_{\lambda} \to K_{\lambda}/N$, we get a map $i \circ \pi \colon K_{\lambda} \to M$ with $M$ being a finite direct sum of $J_{\mu}$. If any of the $J_{\mu}$ occurring in $M$ satisfies $\mu \geq \lambda$, then by composing with the canonical projection $M \to J_{\mu}$ we get a map $K_{\lambda} \to J_{\mu}$ which is not injective, so zero by Corollary~\ref{cor:genericmap}. Therefore, we may assume $M$ is a direct sum of $J_{\mu}$ with $|\mu| < |\lambda|$. Since $N = \ker(i \circ \pi)$ and the object $H_{\lambda}$ is the intersection of the kernels of maps to $I_{\mu}$ with $|\mu| < |\lambda|$, we get $H_{\lambda} = N$, which implies that $H_{\lambda}$ is simple.
\end{proof}

\begin{theorem}\label{thm:genericthm}
\leavevmode
\begin{enumerate}
    \item All objects in $\Mod_S^{\gen}$ are locally of finite length;
    \item the $J_{\lambda}$'s are all the distinct indecomposable injectives in $\Mod_S^{\gen}$;
    \item the objects $H_{\lambda}$ are all the distinct simple objects in $\Mod_S^{\gen}$;
    \item every finite length object in $\Mod_S^{\gen}$ has a finite injective resolution;
    \item the object $H_{\lambda}$ has multiplicity one in $K_{\lambda}$, and every other irreducible component $H_{\mu}$ of $K_{\lambda}$ satisfies $|\mu| < |\lambda|$.
\end{enumerate}
\end{theorem}
\begin{proof}
We proceed by induction on $|\lambda|$ to prove (5) with the claim being trivially true when $|\lambda| < 0$. As in the proof of the previous lemma, we have an injection $K_{\lambda}/H_{\lambda} \to M$ where $M$ is a finite direct sum of $I_{\mu}$ with $|\mu| < |\lambda|$. The object $J_{\mu}$ has a finite filtration by $K_{\nu}$ with $|\nu| = |\mu|$ since the $\GL$-representation $I_{\mu}$ is of finite length (and the category $\Rep^{\pol}(\GL)$ is graded). By induction, all $H_{\nu}$ occurring in $M$ satisfy $|\nu| < |\lambda|$ proving (5) for $\lambda$ as well.

By a similar induction on $\lambda$, we also see that $K_{\lambda}$ (and therefore $J_{\lambda}$) has finite length for all $\lambda$. Every finitely generated object in $\Mod_S^{\gen}$ embeds into a finite direct sum of $J_{\lambda}$ by Proposition~\ref{prop:finalemb}, so every finitely generated object has finite length, proving (1).

Assume $L$ is a simple object in $\Mod_S^{\gen}$. Being finitely generated, it embeds into a finite direct sum of $I_{\lambda}$'s and hence maps into the socle of one of the $I_{\lambda}$ whence $L \cong H_{\lambda}$. Furthermore, for distinct $\lambda$ and $\mu$, the $\GL$-representations $I_{\lambda} \not\simeq I_{\mu}$ and so $J_{\lambda} \not\simeq J_{\mu}$, using which we obtain that $H_{\lambda} \not\simeq H_{\mu}$. This proves (3), and in turn (2) since the indecomposable injectives in a finite length category are precisely the injective envelopes of the simple objects. 

Finally, to prove (4), we first note that the module $K_{\lambda}$ has finite injective dimension for each $\lambda$. Indeed, we may start with a finite injective resolution $0 \to I_0 \to I_1 \to \ldots \to I_n \to 0$ of the simple $\GL$-representation $L_{\lambda}$ and apply the exact functor $T(S \otimes -)$ to obtain a finite injective resolution of $K_{\lambda}$. We proceed by induction on $|\lambda|$ again to show that $H_{\lambda}$ has a finite injective resolution. Given an $H_{\lambda}$, we have an exact sequence
\[ 0 \to H_{\lambda} \to K_{\lambda} \to K_{\lambda}/H_{\lambda} \to 0.\]
The two objects $K_{\lambda}$ and $K_{\lambda}/H_{\lambda}$ have finite injective resolutions -- the former by the previous paragraph, and the latter by using (5) and induction, so $H_{\lambda}$ also has a finite injective resolution. By (3), we see that all finite length objects have a finite injective resolution, as required.
\end{proof}

\section{Shift Theorem}\label{s:nagpalshift}
We prove the shift theorem (Theorem~\ref{thm:shiftintro}) at the end of this section. We assume $q = p^r$ with $r > 0$ throughout.

\subsection{Hasse--Schur derivatives}\label{ss:hasseschur} 
In this section, we define a sequence of endofunctors of $\Rep^\pol(\GL)$ called the \textit{Hasse--Schur derivatives}. We first define it for the category of strict polynomial functors $\Pol$ and then define it for $\Rep^\pol(\GL)$ using the equivalence $\Phi \colon \Rep^{\pol}(\GL) \to \Pol$ discussed in \cite[Section~2.1]{gan22ext}. These functors restrict to give endofunctors of $\Mod_A$ for any $\GL$-algebra $A$.

\subsubsection{Hasse--Schur derivatives of polynomial functors} Fix $m \geq 0$. Let $F \colon \Fec_k \to \Fec_k$ be a polynomial functor. We define the \textit{$m$-th Hasse--Schur derivative} of $F$, denoted $\Sh_m(F)$, to be the polynomial functor which on objects is defined by
\begin{displaymath}
\Sh_m(F)\{V\} = F\{k \oplus V\}^{[m]},
\end{displaymath}
where the superscript here denotes the subspace of $k \oplus V$ on which the top one-dimensional torus $\bG_m$ acts with weight $m$. Given a linear map $f \colon V \to W$, the induced map $\Sh_m(F)(f)$ is given by restricting the map $F(\id_k \oplus f) \colon F\{k \oplus V\} \to F\{k \oplus W\}$ to the appropriate subspace. It is easy to see that $\Sh_m(F)$ is a polynomial functor and that the assignment $F \to \Sh_m(F)$ is functorial. 

\subsubsection{Basic properties of the Hasse--Schur derivative}
For all $m$, the $m$-th Hasse--Schur derivative is an exact functor (because $\bG_m$ is semisimple). It preserves finite length objects of $\Pol$. Furthermore, the sequence $\{\Sh_m\}$ is a categorification of the usual Hasse derivatives (of polynomials) as 
\begin{itemize} 
    \item $\Sh_0$ is the identity functor,
    \item for all $m$, the functor $\Sh_m$ is a $k$-linear endofunctor, and 
    \item the sequence satisfies the generalized Leibniz rule, i.e., $\Sh_m(F\otimes G) \cong \bigoplus_{i+j = m} \Sh_i(F) \otimes \Sh_j(G)$.
\end{itemize}
Furthermore, since all the weights occurring in an $r$-fold Frobenius twisted representation are divisible by $p^r$, we also see that $\Sh_m(W^{(r)}) = 0$ for $0 < m < p^r$.

\subsubsection{Hasse--Schur derivatives of a polynomial representation}Since $\Pol$ and $\Rep^\pol(\GL)$ are equivalent abelian categories, we obtain a sequence of endofunctors of $\Rep^\pol(\GL)$, which we again denote by $\{\Sh_m \}$. For a polynomial representation $W$, we identify $\Sh_m(W)$ as a subspace of $W\{k \oplus \bV\}$, where $\bG_m$ acts with weight $m$. The newly introduced basis vector of $k \oplus \bV$ will usually be denoted by $f$ or $y$ (or possibly with subscripts on these variables) depending on the situation. 

\subsubsection{Hasse--Schur derivative of an \texorpdfstring{$A$}{A}-module} Let $A$ be a $\GL$-algebra, and $M$ be an $A$-module. The Hasse--Schur derivative $\Sh_m(M)$ is canonically an $A$-module, as we now explain. We have a canonical $\GL$-equivariant map $A \to \Phi(A)(k \oplus \bV)$ induced by the canonical inclusion $\bV \to k \oplus \bV$. The $\GL$-representation $\Phi(M)(k \oplus \bV)$ is a $\GL(k \oplus \bV)$-equivariant module over $\Phi(A)(k \oplus \bV)$, so by restriction of scalars, it is a $\GL$-equivariant module over $A$. The subspace $\Sh_m(M)$ of $\Phi(M)(k \oplus V)$ is stable under the action of $\GL$ as well as multiplication by elements of $A$, so $\Sh(M)$ is an $A$-module.

To summarize: for all $m \geq 0$, the $m$-th Hasse--Schur derivative induces a functor 
\[\Sh_m \colon \Mod_A \to \Mod_A.\]

\subsection{The shift functor}\label{ss:shiftdefn}
For the remainder of this section, we work with $\Mod_{\Sm}$. 
We now define a natural map $i_M \colon M \to \mShq(M)$ for an $\Sm$-module $M$. We have an $\Sm$-module map 
\begin{displaymath}
m \mapsto y_1^{q/p} m
\end{displaymath}
from $M$ to $M\{k \oplus \bV\}$, where $y_1$ is the new variable in $\Sm\{k \oplus \bV\}$. Since $m$ uses only the basis vectors $e_i$, the top $\bG_m \subset \GL(k \oplus \bV)$-weight of $y_1^{q/p}m$ is $q/p$. So the map factors $M \to \mShq(M) \to M\{k^r \oplus \bV\}$. We let $i_M$ be the first map. It is easy to see that $i$ is natural. 
We let $\mDeq(M) = \coker(i_M)$ be the \textit{difference functor}. Since $\mShq$ is an exact functor, it follows from the snake lemma that $\mDeq$ is right exact.

We first prove some basic properties of the Hasse-Schur derivative functor (which hereafter we refer to as the \textit{shift functor}). We also let $\mGaq$ be the functor that takes an $\Sm$-module to its torsion submodule.
\begin{proposition}\label{prop:basicshift}
    Let $M$ be an $\Sm$-module.
    \begin{enumerate}[label=(\alph*)]
        \item The kernel of $i_M$ is contained in the torsion submodule of $M$.
        \item The map $i_M$ is injective if and only if $M$ is torsion-free.
        \item If $M$ is torsion-free, then so is $\mShq(M)$.
        \item If $M$ is finitely generated, then $\mShq^r(M)$ is torsion-free for sufficiently large $r$.
        \item If $M$ is nonzero, the module $\mDeq(M)$ is generated in degrees $\le t_0(M) - 1$.
        \item If $M$ is a finitely generated semi-induced module, then so are $\mShq(M)$ and $\mDeq(M)$.
        \item If $M$ is finitely generated, then so are $\mShq(M)$ and $\mDeq(M)$.
        \item The functors $\mShq\mDeq$ and $\mDeq\mShq$ are naturally isomorphic (i.e., the functors $\mShq$ and $\mDeq$ commute).
        \item The functors $\mShq \mGaq$ and $\mGaq\mShq$ are naturally isomorphic (i.e., the functors $\mShq$ and $\mGaq$ commute).
    \end{enumerate}
\end{proposition}

The proposition will be proved at the end of this section, after we prove a few preliminary results. 

\begin{lemma}\label{lem:obvshiftlemma}
    Let $i \colon \id \to \mShq$ be the natural map of functors. We have two natural maps from $\mShq$ to $\mShq^2$: one given by $\mShq (i)$, and the other given by $i_{\mShq}$. There exists an involution $\tau$ of $\mShq^2$ such that $\mShq(i)=\tau \circ i_{\mShq}$.
\end{lemma}
\begin{proof}
The proof of \cite[Lemma~4.2]{gan22ext} applies mutatis mutandis.
\end{proof}

\begin{lemma} \label{lem:shiftofinduced}
    Let $V$ be a polynomial representation of $\GL$. We have isomorphisms
   \begin{displaymath}
   \mShq(\Sm \otimes V) \cong \bigoplus_{i=0}^{{q}/{p}} \Sm \otimes \Sh_m(V).
   \end{displaymath}
   Furthermore, the natural map $i_{\Sm \otimes V}$ is the identity map onto the first summand. Therefore, we have isomorphisms
  \begin{displaymath}
   \mDeq(\Sm\otimes V) \cong \bigoplus_{i=1}^{{q}/{p}} \Sm \otimes \Sh_m(V).
  \end{displaymath} 
\end{lemma}
\begin{proof}
We have isomorphisms $\Sh_m(\Sm) \cong \Sm$ for $m < q$. We obtain the result using the generalized Leibniz rule (see also the proof of \cite[Lemma~4.3]{gan22ext}).
\end{proof}
\begin{corollary}\label{cor:dfrobzero}
    Let $W$ be a polynomial representation of $\GL$. Then $\mDeq(\Sm\otimes W^{(r)}) = 0$.
\end{corollary}
\begin{proof}
All the weights of $W^{(r)}$ are divisible by $p^r$, so $\Sh_m(W^{(r)}) = 0$ for all $m < q$ . The result now follows from Lemma~\ref{lem:shiftofinduced}.
\end{proof}
\begin{proof}[Proof of Proposition~\ref{prop:basicshift}]
We will prove parts a, b and d; for the remaining, the proof of a similar proposition \cite[Proposition~4.1]{gan22ext} applies.

a) Assume $y_1^{q/p}m = 0$ in $\mShq(M)$ for some nonzero $m \in M$. As the module $\mShq(M)$ is contained in $M\{k \oplus \bV\}$, we have $y_1^{q/p}m = 0$ in $M\{k \oplus \bV\}$, which is a $\GL(k \oplus \bV)$-equivariant module over $\Sm\{k \oplus \bV\}$. Since $y_1^{q/p}$ and $m$ are disjoint, we get that $\langle y_1^{q/p} \rangle \langle m \rangle = 0$ by Lemma~\ref{lem:disjointnzd}. In particular, we get that the ideal $\langle x_1^{q/p} \rangle \subset \Sm$ annihilates $m \in M$, and so $m$ is torsion.

b) If $M$ is torsion-free, then by the previous part, we have $\ker(i_M)=0$. Now assume $M$ is not torsion-free. By Lemma~\ref{lem:locannihiltor}, we may choose a nonzero element $m$ annihilated by $\fm^{q/p}$. Choose $r \gg 0$ such that $x_r$ and $m$ are disjoint. Then $x_r^{q/p} m = 0$. The group $\GL(k \bigoplus \bV)$ acts on $\Sm\{k \oplus \bV\}$, so we may apply the element $g$ of $\GL(k \oplus \bV)$ which sends $e_r \to e_r + f_1$ and fixes every other basis vector to the equation $x_r^{q/p} m = 0$ to obtain that $y_1^{q/p} m = 0$, which implies that $m \in \ker(i_M)$. 

d) First, assume $M$ is a finitely generated $S/\fm^{[q/p]}$-module of degree $\leq d$, so that it is a quotient of $S/\fm^{[q/p]} \otimes V$ with $V$ a $\GL$-representation of degree $\leq d$. Since $S/\fm^{[q/p]}$ has no weight vectors with a component $= q/p$, we get $\mShq(S/\fm^{[q/p]}) = 0$, and combining this with the generalized Leibniz rule, we get $\mShq(S/\fm^{[q/p]} \otimes V)$ is of the form $S/\fm^{[q/p]} \otimes W$ with $W$ a $\GL$-representation of degree $\leq d-1$. By induction, we see $\mShq^r(M)$ vanishes for $r > d$. Now, given a finitely generated torsion $\Sm$-module $M$, using Lemma~\ref{lem:locannihiltor}, we get that $M$ has a finite filtration by $S/\fm^{[q/p]}$-modules, and since $\mShq$ is an exact functor, the module $\mShq^r(M)$ vanishes for $r \gg 0$ as claimed.
\end{proof}

\begin{lemma}\label{lem:handlingtorsion}
{   Let $0 \to L \to M \to N \to 0$ be a short exact sequence of $\Sm$-modules. We have a six-term exact sequence
   \begin{displaymath}
       0 \to \bKq(L) \to \bKq(M) \to \bKq(N) \to \mDeq(L) \to \mDeq(M) \to \mDeq(N) \to 0
   \end{displaymath}
   where $\bKq = \ker(\id \to \mShq)$.}
\end{lemma}
\begin{proof}
{
    The result follows by applying the snake lemma to the diagram}
\begin{displaymath}
\begin{tikzcd}
  0 \arrow[r] & L \arrow[d, "i_L"] \arrow[r] & M \arrow[d, "i_M"] \arrow[r] & N\arrow[d, "i_N"] \arrow[r] & 0 \\
  0 \arrow[r] & \mShq(L) \arrow[r] & \mShq(M) \arrow[r] & \mShq(N) \ar[r] & 0.
\end{tikzcd}
\end{displaymath}
\end{proof}

\begin{corollary}\label{cor:handlingtorsion}
Let $0 \to L \to M \to N \to 0$ be a short exact sequence of $\Sm$-modules such that $N$ is torsion-free. We have a short exact sequence 
\begin{displaymath}
0 \to \mDeq(L) \to \mDeq(M) \to \mDeq(N) \to 0
\end{displaymath}
\end{corollary}
\begin{proof}
The module $\bKq(N) = 0$ by Proposition~\ref{prop:basicshift}(a). The result follows by using the previous lemma.
\end{proof}
\subsection{Vanishing of \texorpdfstring{$\mDeq$}{Delta}} \label{ss:deltavanishing}
In this section, we prove a partial converse to Corollary~\ref{cor:dfrobzero}: we show that if $M$ is a torsion-free $\Sm$-module with $\mDeq(M) = 0$, then it must be an extension of $\Sm \otimes W^{(r)}$ with $W$ varying. In particular, such a module must be semi-induced. The strategy of proof here is the one we employ for modules over the infinite exterior algebra \cite[Section~4.3]{gan22ext} so a reader who finds this section technical might find reading loc.~cit.~first more clarifying.

The results of this section hinge on the next lemma. The proof idea is that we can find an element $v$ in degree $n$ of weight $\lambda$ with $\lambda_1 < q$. The element $y_1^{q/p - \lambda_1} \sigma(v)$ is an element in $\mDeq(M)$ where $\sigma$ is a permutation that swaps the newly added basis vector $f$ with $e_1$. The hard part is to show that the element $y_1^{q/p - \lambda_1}\sigma(v)$ is nonzero, which is where the assumption of torsion-free comes into play.
\begin{lemma}\label{lem:dmzero}
    Let $M$ be a torsion-free finitely generated $\Sm$-module with $M_i = 0$ for $i < n$ and $M_n$ not $r$-fold Frobenius twisted (i.e., $M_n$ is not of the form $U^{(r)}$ for a polynomial representation $U$). Then $\mDeq(M) \ne 0$. 
\end{lemma}
\begin{proof}
    It suffices to exhibit a nonzero element of degree $< n$ in $\mShq(M)$ since the image of the natural map $M \to \mShq(M)$ lies in degree $\geq n$. So, assume that the $\GL$-representation $M_n$ is not $r$-fold Frobenius twisted. Choose a minimally flat weight vector $v \in M_n$ (see Definition~\ref{def:minflatvec}). By using a suitable permutation, we may assume that $v$ has weight $\lambda$ with $\lambda_1 = p^i < q$ and every nonzero component $\lambda_j$ of $\lambda$ satisfies $\lambda_j \geq p^i$. Since $v$ is minimally flat, given a weight $\mu$ occurring in $M_n$, the nonzero components of $\mu$ are $\geq p^i$.
    
    In $M\{k \oplus \bV\}$, let $w = \sigma v$ where $\sigma$ is a permutation that swaps the basis vectors $f$ and $e_1$ and fixes the other $e_i$. 
    Consider the element $y_1^{q/p - p^i} w \in M\{k \oplus \bV\}$. The $\bG_m$ weight of $y_1^{q/p-p^i}w$ is clearly $q/p$ so it lies in $\mShq(M)$ and has degree $ n - p^i < n$. So we will be done if we can show $y_1^{q/p-p^i}w$ is nonzero, which is trivial if $i = r - 1$ or $p^i = q/p$. So assume $0 \leq i < r - 1$, and for the sake of contradiction, that $y_1^{q/p-p^i} w = 0$. 
    
    Let $W$ be the 2-dimensional subspace of $k \oplus \bV$ spanned by $f$ and $e_1$. We identify the subgroup $GL(W)$ with $\GL_2$; the subgroup acts on $M\{k \oplus \bV\}$, and therefore, so does the algebra of distributions of $\GL_2$. Let $M'$ be the $\GL_2$ representation generated by $w$. We let $E = e_{1,2}^{(p^i)} \in \Dist(\GL_2)$. By applying $E$ to the equation $y_1^{q/p-p^i}w = 0$, we obtain the equation
    \begin{displaymath}
        \sum_{j=0}^{p^i} e_{1,2}^{(j)}(y_1^{q/p-p^i}) e_{1,2}^{(p^i - j)}(w) = 0.
    \end{displaymath} 
   All the terms $e_{1, 2}^{(j)}(y_1^{q/p-p^i})$ with $0 < j < p^i$ vanish by Lemma~\ref{lem:hassevanishing}, so the above equation simplifies to $E(y_1^{q/p-p^i}) w + y_1^{q/p-p^i}E(w) = 0$. By multiplying this equation with $y_1^{q/p-p^i}$ and using our assumption that $y_1^{q/p-p^i}w = 0$, we also get that 
   \[y_1^{2(q/p-p^i)} E(w) = 0.\]
   Since $2(q/p - p^i) < q$, the element $y_1^{2(q/p - p^i)} \ne 0$ in $\Sm\{k \oplus \bV\}$. The $\GL_2$-representation $M'$ is generated by a weight vector of weight $(0, p^i)$, and by assumption on $p^i$, all the nonzero components of the weights occurring in $M$, and in turn $M'$ are $\geq p^i$. Therefore, the only weights occurring in $M'$ are of the form $(0, p^i)$ or $(p^i,0)$ which implies that the $\GL_2$-representation $M' \cong ((k^2)^{(i)})^{\oplus t}$ for some $t \in \bN$. So $E(w)$ is nonzero and it has weight $(0, p^i)$ for the $\GL_2$-action. So we obtain two nonzero disjoint elements $y_1^{2(q/p-p^i)}$ and $E(w)$ with $y_1^{2(q/p-p^i)}E(w) = 0$ which implies that $E(w)$ is a torsion element by Lemma~\ref{lem:disjointnzd} but this contradicts the assumption that $M$ is torsion-free.
\end{proof}
The next result is similar to \cite[Proposition~4.10]{gan22ext}.
\begin{proposition}\label{prop:dmzerogenonedeg}
    Let $M$ be a torsion-free $\Sm$-module generated in degree $n$ such that $\mDeq(M) = 0$. The natural map $\Sm \otimes M_n \to M$ is an isomorphism, i.e., $M$ is an induced $\Sm$-module.
\end{proposition}
\begin{proof}
    The irreducible components of $\Sm$ in positive degree are $p^r$-restricted representations of $\GL$ as all the weights are clearly $p^r$-restricted. Therefore, if $W$ is an irreducible representation of $\GL$, then $(\Sm)_n \otimes W^{(r)}$ is not $r$-fold Frobenius twisted for $n>0$. It follows that for a finite length polynomial representation $W$, the irreducible components of $\Sm \otimes W^{(r)}$ are not $r$-fold Frobenius twisted in degrees $ > \deg(W^{(r)})$.
    
    By the assumptions on $M$, we see that $M_n$ is $r$-fold Frobenius twisted by Lemma~\ref{lem:dmzero}. Consider the natural map $\phi \colon \Sm \otimes M_n \to M$. The map $\phi$ is surjective, and $K = \ker(\phi)$ is zero in degrees $\le n$. We have to prove that $K = 0$. First, we claim that $\mDeq(K) = 0$. Indeed, since $M$ is torsion-free, by Corollary~\ref{cor:handlingtorsion}, we have a short exact sequence 
     \begin{displaymath}
    0 \to \mDeq(K) \to \mDeq(\Sm \otimes M_n) \to \mDeq(M) \to 0.
    \end{displaymath} 
    By Corollary~\ref{cor:dfrobzero}, we have $\mDeq(\Sm \otimes M_n) = 0$ since $M_n$ is $r$-fold Frobenius twisted, and so $\mDeq(K) = 0$, as claimed. 
    Now, suppose $K$ is nonzero, and let $s$ be the smallest degree such that $K_s \ne 0$. Since $s > n$ the $\GL$-representation $K_s$ is not $r$-fold Frobenius twisted by the previous paragraph. Therefore, by Lemma~\ref{lem:dmzero}, we see that $\mDeq(K)$ must be nonzero, which is a contradiction. Therefore, the map $\phi$ is an isomorphism, as required.
\end{proof}

\begin{lemma}\label{lem:torsionfreequotient}
    Let $M$ be a torsion-free $\Sm$-module with $t_0(M) = n$ such that $\mDeq(M)$ is semi-induced and $\mDeq(M/M^{<n}) = 0$. Then $M/M^{<n}$ is torsion-free.
\end{lemma}
\begin{proof}
    By Corollary~\ref{cor:handlingtorsion}, we have the short exact sequence 
    \begin{displaymath}
    0 \to \mKq(M/M^{<n}) \to \mDeq(M^{<n}) \to \mDeq(M) \to 0.
    \end{displaymath} 
    Applying $\Tor$ to this gives us the exact sequence
    \begin{displaymath}
    \Tor_1^{\Sm}(\mDeq(M), k) \to \Tor_0^{\Sm}(\mKq(M/M^{<n}),k) \to \Tor_0^{\Sm}(M^{<n},k).
    \end{displaymath}
    However, since $\mDeq(M)$ is semi-induced, by Proposition~\ref{prop:flatequalssemi}, we get $\Tor_1^{\Sm}(\mDeq(M), k) = 0$, which implies that 
    \begin{displaymath}
    t_0(\mKq(M/M^{<n})) \leq t_0(M^{<n}) < n.
    \end{displaymath}
    But $\mKq(M/M^{<n})$ is supported only in degrees $\ge n$ (as it is a submodule of $M/M^{<n}$), and so must be zero. Therefore, the natural map $M/M^{<n} \to \mShq(M/M^{<n})$ is injective, which implies that the module $M/M^{<n}$ is torsion-free by Proposition~\ref{prop:basicshift}(b).
\end{proof}
We can now prove the result promised at the beginning of this section.
\begin{proposition}\label{prop:dmzeroisflat}
Let $M$ be a finitely generated torsion-free $\Sm$-module such that $\mDeq(M) = 0$. Then $M$ is semi-induced.
\end{proposition}
\begin{proof}
    We proceed by induction on the generation degree of $M$.
    When $t_0(M) = 0$, this follows from Proposition~\ref{prop:dmzerogenonedeg}.
    Assume that $t_0(M) = n > 0$. Since $\mDeq$ is right exact, we see that $\mDeq(M/M^{<n}) = 0$. Therefore, by Lemma~\ref{lem:torsionfreequotient}, the module $M/M^{<n}$ is torsion-free, and so by Proposition~\ref{prop:dmzerogenonedeg}, the module $M/M^{<n}$ is semi-induced. By Corollary~\ref{cor:handlingtorsion}, we have a short exact sequence 
    \begin{displaymath}
    0 \to \mDeq(M^{<n}) \to \mDeq(M) \to \mDeq(M/M^{<n}) \to 0.
    \end{displaymath}
    As $\mDeq(M)=0$, we have that $\mDeq(M^{<n}) = 0$, which implies that $M^{<n}$ is semi-induced (by induction) and therefore, so is $M$. 
\end{proof}
We also note a corollary of Lemma~\ref{lem:torsionfreequotient} that we will use in the proof of the shift theorem.
\begin{corollary}
 \label{cor:dqlowdeg}
     Let $M$ be a finitely generated torsion-free $\Sm$-module with $t_0(M) = n$ such that $\mDeq(M)$ is semi-induced with $t_0(\mDeq(M)) < \max(-1, n - \frac{q}{p})$. Then $M/M^{<n}$ is semi-induced. 
 \end{corollary}
 \begin{proof}
     Since the functor $\mDeq$ is right exact, we have \begin{displaymath}
     t_0(\mDeq(M/M^{<n})) \le t_0(\mDeq(M)) < n-\frac{q}{p}.
     \end{displaymath}
     Since $M/M^{<n}$ is supported in degrees $\ge n$, the module $\mShq(M/M^{<n})$ is supported in degrees $\ge n-\frac{q}{p}$, and so $\mDeq(M/M^{<n})$ is also supported in degrees $\geq n-\frac{q}{p}$. Therefore, the above inequality implies that $\mDeq(M/M^{<n}) = 0$. By Lemma~\ref{lem:torsionfreequotient}, the module $M/M^{<n}$ is torsion-free and by applying Proposition~\ref{prop:dmzerogenonedeg}, we see that the module $M/M^{<n}$ is semi-induced.
 \end{proof}

\subsection{Semi-induced subquotients of semi-induced modules}\label{ss:subofsemi}
We prove some lemmas bounding the generation degree of semi-induced subquotients of semi-induced modules. The only idea used in this section is that $\Sm\{k^n\}$ is supported in degrees $0$ through $n(q-1)$. For an $\bN$-graded object $M$, we let $\maxdeg(M)$ be the largest integer such that $M_i \ne 0$, with $\maxdeg(0) = -1$ and $\maxdeg(M) = \infty$ if no such integer exists.

\begin{lemma}
	Let $M$ be an $\Sm$-module. Then for all $n$, $\maxdeg(M\{k^n\}) < (q-1)n + t_0(M)$. Furthermore, if $M$ is semi-induced, then $\maxdeg(M\{k^n\}) = {(q-1)n + t_0(M)}$ for all sufficiently large $n$.
\end{lemma}
\begin{proof}
    This is clear since $\Sm$ evaluated at $k^n$ is supported in degrees $\le (q-1)n$.
\end{proof}

\begin{corollary}\label{cor:subsemi}
Let $F$ be a semi-induced module and $Z$ be a submodule of $F$ such that $Z/Z^{< t_0(Z)}$ is semi-induced. Then $t_0(Z) \leq t_0(F)$.
\end{corollary}
\begin{proof}
    For sufficiently large $n$, $\maxdeg(Z\{k^n\}) = (q-1)n + t_0(Z)$ as $Z/Z^{< t_0(Z)}\{k^n\}$ is nonzero in that degree by the previous lemma. Since $Z$ is a submodule of $F$, we have $\maxdeg(Z\{k^n\}) \le \maxdeg(F\{k^n\})$ for sufficiently large $n$, or $(q-1)n + t_0(Z) \leq (q-1)n + t_0(F)$ for large $n$, giving us the required inequality.
\end{proof} 

\begin{lemma}
    Let $0 \rightarrow W \rightarrow F \rightarrow M \rightarrow 0$ be an exact sequence of $\Sm$-modules with $F$ semi-induced and $t_0(F) = t_0(M)$. If $t_0(W) \leq t_0(F)$, then $M/M^{<t_0(M)}$ is semi-induced.
\end{lemma}
\begin{proof}
Let $n = t_0(M)$. It is easy to see that since $t_0(W) \le t_0(F)$, we have $t_1(M) \leq t_0(M)$. Now consider the exact sequence,
\begin{displaymath}
    0 \to M^{<n} \to M \to M/M^{<n} \to 0.
\end{displaymath}
The long exact sequence corresponding to $\Tor$ gives us
\begin{displaymath}
 \Tor_1^{\Sm}(M,k)  \to  \Tor_1^{\Sm}(M/M^{<n},k) \to \Tor_0^{\Sm}(M^{<n}, k) 
\end{displaymath}
which implies that $t_1(M/M^{<n}) \le \max(t_1(M), t_0(M^{<n}) \le t_0(M) = n$. Therefore, by Lemma~\ref{lem:relationsinlowdeg}, we have that $M/M^{<n}$ is semi-induced, as required.
\end{proof}

\begin{corollary}\label{cor:descentinequality}
    Let $M$ be an $\Sm$-module with a degree-minimal flat resolution
    \begin{displaymath}
        \ldots \to F_j \to F_{j-1} \to \ldots F_1 \to F_0 \to M \to 0
    \end{displaymath}
    such that $t_0(F_{i+1}) \leq t_0(F_{i})$ for some $i \ge 0$. Then $t_1(M) \leq t_0(M)$.
\end{corollary}
\begin{proof}
     Let $Z_{j+1}$ be the kernel of the map from $F_j$ to $F_{j-1}$ (with $F_{-1} = M = Z_0$). We have the short exact sequence
    \begin{displaymath}
        0 \to Z_{i+1} \to F_i \to Z_{i} \to 0
    \end{displaymath}
    with $t_0(F_i) = t_0(Z_i)$.
    By assumption, we have $t_0(Z_{i+1}) = t_0(F_{i+1}) \leq t_0(F_i)$ and so by the previous lemma $Z_i/(Z_i^{< t_0(Z_i)})$ is semi-induced. The module $Z_i$ is a submodule of the semi-induced module $F_{i-1}$ with $Z_i/(Z_i^{< t_0(Z_i)})$ semi-induced. So using Corollary~\ref{cor:subsemi}, we obtain the inequality $t_0(Z_i) \leq t_0(F_{i-1})$, and since $t_0(Z_i) \le t_0(F_i)$, we get that $t_0(F_i) \leq t_0(F_{i-1})$ as well. By iterating this, we get that $t_0(F_{j+1}) \leq t_0(F_j)$ for all $0 \leq j \leq i$ and in particular, that $t_0(F_1) \le t_0(F_0)$. By Proposition~\ref{prop:degreeres}, $t_1(M) \leq \max(t_1(M), t_0(M)) = \max(t_0(F_1), t_0(F_0)) = t_0(F_0) = t_0(M)$, as required. 
\end{proof}

\subsection{Proof of the shift theorem}\label{ss:pfshift}
We note a result about free resolutions over complete intersections that will be used in our proof.

\begin{proposition}\label{prop:cibound}
    Fix integers $d, r \geq 1$ and let $I$ be an ideal generated by a regular sequence of $r$ homogeneous polynomials of the same degree $d$ in a polynomial ring $Q = k[x_1,x_2,x_3, \ldots, x_r]$. Let $M$ be a finitely generated graded module over the artinian ring $Q/I$ and set $t_i(M) = \maxdeg \Tor_i^{Q/I}(M,k)$. Assume further that $t_i(M) < bi + c$ for real numbers $b, c$ with $b < \frac{d}{2}$. Then $M$ is flat.
\end{proposition}
\begin{proof}
Assume $M$ is not flat. Then $M$ has infinite projective dimension by the Auslander--Buchsbaum formula. We may replace the module $M$ by a sufficiently large syzygy of $M$ as the assumptions continue to hold for this syzygy module. Eisenbud--Peeva \cite[Corollary~6.2.1]{ep16mfrci} have computed the graded Poincar{\'e} series of such a large syzygy (in loc.~cit.~the coefficient of $x^iz^j$ is $\dim_k(\Tor_i(M,k)_j)$.) There exist polynomials $\{m_{j,0}(z)\}_{0 \le j \le r}$ with non-negative integer coefficients that are independent of $i$ such that the coefficient of $x^{2i}$ is 
\begin{displaymath}
    z^{di} [m_{r, 0}(z) + \binom{i + 1}{1} m_{r-1, 0}(z) + \binom{i+2}{2} m_{r-2, 0} + ... + \binom{i+r-1}{r-1} m_{1, 0}(z)].
    \end{displaymath}
    Let $N = \max_j \deg m_{j, 0}(z)$. For $i > n$, the polynomial 
       $m_{n, 0}(z) + \binom{i + 1}{1} m_{n-1, 0}(z) + \binom{i+2}{2} m_{n-2, 0} + ... + \binom{i+n-1}{n-1} m_{1, 0}(z)$ has degree $N$ as all the polynomials have non-negative coefficients and so for $i > n$, the coefficient of $x^{2i}$ is a polynomial in $z$ of degree $di + N$. In particular, for $i \gg 0$, we have $t_{2i}(N) = di + N$ which contradicts the fact $b < \frac{d}{2}$.
\end{proof}

The next lemma puts together what we have proved in the previous two sections. We refer the reader to Section~\ref{sss:shiftintro} for some motivation regarding this lemma.

\begin{lemma}\label{lem:t1lesst0}
    Let $M$ be a torsion-free $\Sm$-module such that $\mDeq(M)$ is semi-induced. Then $t_1(M) \leq t_0(M)$.
\end{lemma}
\begin{proof}
   Let 
   \begin{displaymath}
       \ldots \to F_i \to F_{i-1} \to \ldots \to F_1 \to F_0 \to 0
   \end{displaymath}
    be a degree-minimal flat resolution of $M$, and let $Z_0 = M$ and $Z_{j+1} = \ker(F_j \to F_{j-1})$ for $j>0$. We have short exact sequences 
    \begin{displaymath}
    0 \rightarrow Z_{j+1} \rightarrow F_j \rightarrow Z_{j} \rightarrow 0.
    \end{displaymath}
    Since all the modules are torsion-free, applying the functor $\mDeq$ to these short exact sequences preserves exactness by Corollary~\ref{cor:handlingtorsion}. So we obtain short exact sequences
    \begin{displaymath}
    0 \rightarrow \mDeq(Z_{j+1}) \rightarrow \mDeq(F_j) \rightarrow \mDeq(Z_j) \rightarrow 0
    \end{displaymath}
    for all $j \geq 0$.
    The module $\mDeq(Z_{0}) = \mDeq(M)$ is semi-induced, and $\mDeq(F_j)$ is semi-induced for all $j$ by Proposition~\ref{prop:basicshift}(f). By induction we see that $\mDeq(Z_j)$ is semi-induced for all $j \geq 0$ by Corollary~\ref{cor:semises}. Using the long exact sequence of $\Tor$, we get,
    \begin{displaymath} 
        t_0(\mDeq(Z_{j+1})) \le t_0(\mDeq(F_j)) 
    \end{displaymath}
    for all $j \geq 1$.
    
    There are two exhaustive cases, in either of which we will show that the conclusion $t_1(M) \leq t_0(M)$ holds.
    
   The first case is when there exists an $i \geq 1$ such that $t_0(\mDeq(Z_i)) < t_0(Z_i) - \frac{q}{p}$. In this case, by Lemma~\ref{cor:dqlowdeg}, the $\Sm$-module $Z_i/(Z_i)^{< t_0(Z_i)}$ is semi-induced, and therefore $t_0(F_i) = t_0(Z_i) \leq t_0(F_{i-1})$ by Corollary~\ref{cor:subsemi}. It then follows from Corollary~\ref{cor:descentinequality} that $t_1(M) \leq t_0(M)$, as required. 
   
   We move to the second case, when the assumption of the first case does not hold, so for all $i\geq 1$, we have $t_0(\mDeq(Z_i)) \geq t_0(Z_i) - \frac{q}{p}$. Using this we obtain,
    \begin{displaymath}
    t_0(F_i) = t_0(Z_i) \leq t_0(\mDeq(Z_i)) + \frac{q}{p} \leq t_0(\mDeq(F_{i-1})) + \frac{q}{p} \leq t_0(F_{i-1}) + \frac{q}{p} - 1;
    \end{displaymath}
    for the penultimate inequality, we use $t_0(\mDeq(Z_{i})) \leq t_0(\mDeq(F_{i-1}))$ which we showed earlier in this proof, and for the last inequality, we use Proposition~\ref{prop:basicshift}(e). By induction, for all $i \geq 1$, we have
    \begin{displaymath}
    t_0(F_i) \le t_0(F_0) + i(\frac{q}{p} - 1).
    \end{displaymath}
    Let $N > t_1(M) + t_0(M)$ be an integer. Evaluating the degree-minimal flat resolution of $M$ on $k^N$, we obtain a flat resolution of $M\{k^N\}$ over $\Sm\{k^N\}$. Using $t_0(F_i) \leq t_0(F_0) + i(\frac{q}{p} - 1)$, we obtain
    \begin{displaymath}
    t_0(F_i\{k^n\}) \le t_0(F_i) \le t_0(F_0) + j(\frac{q}{p} -1) = t_0(M\{k^n\}) + j(\frac{q}{p} - 1).
    \end{displaymath}
    Since $t_i(M\{k^n\}) \leq t_0(F_i\{k^n\})$, we have, $t_{i}(M\{k^n\}) < t_0(M\{k^n\}) + j(\frac{q}{p} - 1)$. By Proposition~\ref{prop:cibound}, this implies that $M\{k^N\}$ is flat as a $\Sm\{k^N\}$-module, which in particular implies that $\Tor_1^{\Sm}(M, k) = 0$, and so $t_1(M) \le t_0(M)$ in this case as well.
\end{proof}

\begin{proposition}\label{prop:liftingsemi}
	Let $M$ be a finitely generated torsion-free $\Sm$-module such that $\mDeq(M)$ is semi-induced. Then $M$ is also semi-induced.
\end{proposition}
\begin{proof}
	We induct on $t_0(M)$. When $t_0(M) = 0$, the module $\mDeq(M) = 0$, and so $M$ is semi-induced by Proposition~\ref{prop:dmzeroisflat}. Assume now $t_0(M) = n > 0$. 
	We have the short exact sequence
	\begin{displaymath}
	0 \to M^{<n} \to M \to M/M^{<n} \to 0.
	\end{displaymath}
	So it suffices to prove that $M^{<n}$ and $M/M^{<n}$ are semi-induced.
	
	We first show that $M/M^{<n}$ is semi-induced. The long exact sequence of $\Tor_{*}^{\Sm}(-, k)$ associated to the above short exact sequence gives us 
	\begin{displaymath}
    \Tor_1^{\Sm}(M, k) \to \Tor_1^{\Sm}(M/M^{<n}, k) \to \Tor_0^{\Sm}(M^{<n}, k)
    \end{displaymath}
	from which we see that $t_1(M/M^{<n}) \leq \max(t_1(M), t_0(M^{<n})) = \max(t_1(M), n-1) = n$; for the last equality, we use Lemma~\ref{lem:t1lesst0} to get $t_1(M) \leq t_0(M) = n$. By Lemma~\ref{lem:relationsinlowdeg}, we see that $M/M^{<n}$ is a semi-induced module. 
	
    We now proceed to show that $M^{<n}$ is also semi-induced. Since $M/M^{<n}$ is semi-induced, it is torsion-free, and therefore, we have a short exact sequence
    \begin{displaymath}
    0 \to \mDeq(M^{<n}) \to \mDeq(M) \to \mDeq(M/M^{<n}) \to 0
    \end{displaymath}
    by Corollary~\ref{cor:handlingtorsion}. By Proposition~\ref{prop:basicshift}(f), we have that $\mDeq(M/M^{<n})$ is semi-induced, and so from the short exact sequence above, we see that $\mDeq(M^{<n})$ is also semi-induced by Corollary~\ref{cor:semises}. Since $t_0(M^{<n}) < n$, by the induction hypothesis, it now follows that $M^{<n}$ is semi-induced, as required. 
\end{proof}

We can now prove the shift theorem.
\begin{proof}[Proof of Theorem~\ref{thm:shiftintro}]
We follow the proof of \cite[Theorem~3.13]{ly17fi} and proceed by induction on $t_0(M)$. It suffices to prove the result for $\mShq^r(M)$ for $r \gg 0$, and so we may additionally assume $M$ is torsion-free by Proposition~\ref{prop:basicshift}(d). When $t_0(M) = 0$, the module $\mDeq(M)=0$ and so $M$ is semi-induced by Proposition~\ref{prop:liftingsemi}. Assume now $t_0(M) = n > 0$. We have a short exact sequence $0 \to M \to \mShq(M) \to \mDeq(M) \to 0$ with $t_0(\mDeq(M)) < n$. By induction, for $l\gg 0$ the module $\mShq^{l}(\mDeq(M))$ is semi-induced. By Proposition~\ref{prop:basicshift}(h) we have $\mShq^{l}(\mDeq(M)) \cong \mDeq(\mShq^{l}(M))$, and so by Proposition~\ref{prop:liftingsemi} we have that $\mShq^l(M)$ is semi-induced, as required. 
\end{proof}

\section{Proof of Theorem~\ref{thm:gensdbmodintro}}\label{s:structure}\label{ss:proofBE}
We prove two useful and related results for $\Mod_{\Sm}$ and $\Mod_S$, which state that one can resolve every finitely generated module by flat modules ``up to torsion"; for $\Sm$, it is a corollary of Theorem~\ref{thm:shiftintro} and for $S$ it is a corollary of Theorem~\ref{thm:genericthm}.

\begin{proposition}[Resolution Theorem for {$\Sm$}-modules]\label{prop:resolutionSm}
Let $M$ be a finitely generated $\Sm$-module. We have a chain complex of $\Sm$-modules
\begin{displaymath}
0 \to M \to P^0 \to P^1 \to \ldots \to P^m \to 0
\end{displaymath}
satisfying the following properties:
\begin{itemize}
\item each $P^i$ is a finitely generated semi-induced module with $t_0(P^i) \le t_0(M) - i$, and
\item the cohomology of this complex is a torsion $\Sm$-module.
\end{itemize}
{Furthermore, given a map of $\Sm$-modules $f \colon M \to N$, we can choose complexes $M \to P^{\bullet}$ and $N \to Q^{\bullet}$ satisfying the above properties, and a map of complexes $\tilde{f}\colon P^{\bullet} \to Q^{\bullet}$ extending $f$.}
\end{proposition}
\begin{proof}
The proof of \cite[Theorem~6.1]{gan22ext} applies mutatis mutandis.
\end{proof}

\begin{proposition}[Resolution Theorem for $S$]\label{prop:resolutionthmS}
    Let $M$ be a finitely generated $S$-module. We have a chain complex of $S$-modules
    \[ 0 \to M \to P^0 \to P^1 \to \ldots \to P^m \to 0\]
    each $P^i$ is a finitely generated torsion-free injective $S$-module and the cohomology groups are torsion $S$-modules.
\end{proposition}
\begin{proof}
    Given a finitely generated module $M$, the object $T(M)$ has a finite injective resolution 
    \[ 0 \to T(M) \to I^0 \to I^1 \to \ldots \to I^m \to 0 \]
    in $\Mod_S^{\gen}$ by Theorem~\ref{thm:genericthm}(4). 
    Applying $\cS$, we obtain the complex
    \[0 \to \cS(T(M)) \to \cS(I^0) \to \cS(I^1) \to \ldots \to \cS(I^m) \to 0.\] 
    The $S$-module $\cS(I^i)$ is semi-induced because injectives in $\Mod_S^{\gen}$ are direct sums of $T(S \otimes I_{\lambda})$ by Proposition~\ref{prop:injective} and $\cS(T(S \otimes I_{\lambda})) \cong S \otimes I_{\lambda}$ by Corollary~\ref{cor:injsat}. The cohomology of this complex is $\rR^i\cS(T(M))$ which by Proposition~\ref{prop:cleanupprop} is isomorphic to $\rR^{i+1}\Gamma(M)$ which is torsion module. So we may set $P^i = \cS(T(I^i))$ and splice this complex with the unit map $M \to \cS(T(M))$ to obtain the requisite complex.
\end{proof}
We note a very pleasant corollary of the above result for $S$-modules.
\begin{theorem}[Shift theorem for $S$-modules]\label{thm:introshiftS}
Let $M$ be a finitely generated $S$-module. For $m \gg 0$, the $S$-module $\Sh_m(M)$ is flat. 
\end{theorem}
\begin{proof}
    Let $0 \to M \to P^0 \to P^1 \to \ldots \to P^r \to 0$ be a chain complex of $S$-modules satisfying the conclusions of the resolution theorem for $S$. The cohomology of this complex is finitely generated and torsion, so it is supported on $\fm^{[q]}$ for $q = p^r \gg 0$. By applying the exact functor $\Sh_n$ for $n \gg 0$ to this complex, the cohomology vanishes, and so we get an exact complex
    \[ 0 \to \Sh_n(M) \to \Sh_n(P^{0}) \to \Sh_n(P^1) \to \ldots \to \Sh_n(P^r) \to 0.\] 
    The result now follows by Corollary~\ref{cor:semises}.
\end{proof}

\begin{remark}
The above theorem is reminiscent of Nagpal's shift theorem from characteristic zero \cite{nag15fi}, which states that for finitely generated $M$, the $S$-module $\Sh_1^m(M)$ is flat for $m \gg 0$. The functor $\Sh_1^m$ carries an action of $\fS_m$ and in characteristic zero, the functor $\Sh_m$ is naturally isomorphic to the $\fS_m$-invariant part of $\Sh_1^m$ (which is a categorification of the relation $\partial_m = \frac{1}{m!}\partial^m$ between the ordinary and Hasse derivatives over $\bQ$). 
\end{remark}
Before we prove Theorem~\ref{thm:gensdbmodintro}, we define some notions of triangulated categories (that eventually end up being identical in our context). 

Given an abelian category $\cA$, we let $\rD^b(\cA)$ be the bounded derived category of $\cA$ (i.e., chain complexes of objects in $\cA$ with finitely generated cohomology). Given a Serre subcategory $\cB \subset \cA$, we have two notions: one is the derived category $\rD^b(\cB)$ and the second is the $\rD^b_{\cB}(\cA)$ which is the triangulated subcategory of $\rD^b(\cA)$ consisting of chain complexes whose cohomology lies in $\cB$. We have a canonical map $\rD^b(\cB) \to \rD^b_{\cB}(\cA)$. This is not always an equivalence; the next lemma gives sufficient criteria for when it is.
\begin{lemma}\label{lem:derivedlemma}
    Let $\cB \subset \cA$ be a Serre subcategory. Assume $\cB$ and $\cA$ have enough injectives and $\cB$ satisfies property (Inj). The map $\rD^b(\cB) \to \rD^b_{\cB}(\cA)$ is an equivalence.
\end{lemma}
\begin{proof}
    Assume $0 \to B \to A \to A' \to 0$ is a short exact sequence in $\cA$ with $B \in ob(\cB)$. Let $I$ be an injective object in $\cB$ containing $B$. The inclusion $B \to I$ extends to a map $A \to I$ since $I$ is also injective in $\cA$ by property (Inj). Therefore we get a map of short exact sequence
     \[
\begin{tikzcd}
   0 \arrow[r] & B \ar[equal]{d} \arrow[r] & A \ar[dashed]{d} \arrow[r] & A' \ar[dashed]{d} \arrow[r] & 0 \\
   0 \arrow[r] & B \arrow[r] & I \arrow[r] & I/B \arrow[r] & 0
\end{tikzcd}
\]
with $B, I$ and $I/B$ in $\cB$, whence we obtain the lemma by part c(ii) in \cite[1.15 Lemma]{kel99ex}.
\end{proof}

\begin{proposition}\label{prop:triangles}
\leavevmode
\begin{itemize}
    \item Let $M$ be a finitely generated $\Sm$-module. We have a triangle $T \to M \to F \to$ in $\rD^b_{\fgen}(\Mod_{\Sm})$ where $T$ is quasi-isomorphic to a bounded complex of finitely generated torsion $\Sm$-modules, and $F$ is quasi-isomorphic to a bounded complex of semi-induced $\Sm$-modules. 
    \item Let $N$ be a finitely generated $S$-module. We have a triangle $T \to M \to F \to$ in $\rD^b_{\fgen}(\Mod_{S})$ where $T$ is quasi-isomorphic to a bounded complex of finitely generated torsion $S$-modules, and $F$ is quasi-isomorphic to a bounded complex of semi-induced $S$-modules. 
\end{itemize}
\end{proposition}
\begin{proof}
    For the first part, let $M \xrightarrow{f} P^{\bullet}$ be a complex satisfying the conclusion of Proposition~\ref{prop:resolutionSm}. The cone of $f$ is a bounded complex with finitely generated torsion cohomology, so it is quasi-isomorphic to a complex of finitely generated torsion $\Sm$-modules by the previous lemma. Therefore, the distinguished triangle $\cone(f)[-1] \to M \to P^{\bullet} \to$ satisfies the requisite properties.
    
    The second part follows by combining Proposition~\ref{prop:resolutionthmS} with Lemma~\ref{lem:derivedlemma} similar to the previous paragraph.
\end{proof}
Before proving Theorem~\ref{thm:gensdbmodintro}, we prove a related result for the derived category of $\Sm$.
\begin{proposition}\label{prop:gensdbSm}
The bounded derived category $\rD^b_{\fgen}(\Mod_{S/\fm^{[p^r]}})$ is generated by the classes of modules $S/\fm^{[p^e]} \otimes L_{\lambda}$ as $\lambda$ varies over all partitions and $e$ varies over the set $\{0, 1, \ldots, r\}$. 
\end{proposition}
\begin{proof}
    We proceed by induction on $r$ where $q = p^r$; the $r = 0$ case is trivial so assume $r > 0$. 
    The classes of all finitely generated $\Sm$-modules generate $\rD^b_{\fgen}(\Mod_{\Sm})$. Using Proposition~\ref{prop:triangles}, we see that the class of any finitely generated module $M$ is generated by finitely generated torsion and semi-induced $\Sm$-modules. Now, every finitely generated torsion $\Sm$-module has a filtration where the graded pieces are $S/\fm^{[q/p]}$-modules, so they are generated by $S/\fm^{[p^e]} \otimes L_{\lambda}$ with $p^e < q$ by the induction hypothesis.
\end{proof}

\begin{proof}[Proof of Theorem~\ref{thm:gensdbmodintro}]
    Using Proposition~\ref{prop:triangles}, we get that the class of a finitely generated $S$-module is generated by the classes of semi-induced modules and the torsion $S$-modules. By the previous proposition, torsion $S$-modules can be generated by $\Sm \otimes L_{\lambda}$ with $q$ and $\lambda$ allowed to vary. Furthermore, every semi-induced module has a finite filtration where the graded pieces are $S \otimes L_{\lambda}$. So these two classes generate $\rD^b_{\fgen}(\Mod_S)$ as claimed.
\end{proof}

\bibliographystyle{alpha}
\bibliography{bibliography}

\end{document}